\documentclass[a4paper,10pt]{amsart}
\usepackage{amsmath}
\usepackage{amssymb}
\usepackage{color}
\usepackage{a4wide}

\theoremstyle{plain}
\newtheorem{theorem}{Theorem}
\newtheorem{corollary}{Corollary}
\newtheorem{lemma}{Lemma}

\newtheorem{remark}{Remark}

\begin{document}

\title[A posteriori estimates for time-periodic parabolic optimal control problems]
{Functional a posteriori error estimates for time-periodic parabolic optimal control problems}

\author{Ulrich Langer}
\author{Sergey Repin}
\author{Monika Wolfmayr}

\address[U. Langer]{Institute of Computational Mathematics,
Johannes Kepler University Linz,
Altenbergerstra\ss e 69, 4040 Linz, Austria}
\email{ulanger@numa.uni-linz.ac.at}

\address[S. Repin]{V. A. Steklov Institute of Mathematics in St. Petersburg,
Fontanka 27, 191011, St. Petersburg, Russia,
and University of Jyv\"{a}skyl\"{a}, Finland}
\email{repin@pdmi.ras.ru}

\address[M. Wolfmayr]{Johann Radon Institute for Computational and Applied Mathematics,
Altenbergerstra\ss e 69, 4040 Linz, Austria}
\email{monika.wolfmayr@ricam.oeaw.ac.at}

\begin{abstract}
This paper is devoted to the {\it a posteriori} error analysis of 
multiharmonic finite element approximations to distributed 
optimal control problems with time-periodic state equations of parabolic type.
We derive {\it a posteriori} estimates of functional type, which are 
easily computable and provide guaranteed upper bounds for
the state and co-state errors as well as for the cost functional.
These theoretical results are confirmed by several numerical tests that
show high efficiency of the {\it a posteriori} error bounds.


\end{abstract}

\maketitle

\section{Introduction}
\label{Sec1:Intro}

We consider optimal control problems with time-periodic pa\-ra\-bolic state equations. 
Problems of this type often arise in different practical applications, e.g., 
in electromagnetics, chemistry, biology, or heat transfer, see also 
\cite{LRW:Altmann:2013, LRW:AltmannStingelinTroeltzsch:2014,LRW:HouskaLogistImpeDiehl:2009}.
Moreover, optimal control problems are the subject matter of lots of different works, see, e.g.,
\cite{LRW:NeittaanmaekiSprekelsTiba:2006, LRW:HinzePinnauUlbrichUlbrich:2009, 
LRW:Troeltzsch:2010, LRW:BorziSchulz:2012}
and the references therein.
The multiharmonic finite element method (MhFEM) is well adapted to the
class of parabolic time-periodic problems.
Within the framework of this method, the state functions are expanded
into Fourier series in time with coefficients depending on the spatial variables. In numerical computations,
these series are truncated and the Fourier coefficients are approximated by the
finite element method (FEM).
This scheme leads to the MhFEM (also called harmonic-balanced FEM), which
was successfully used for the simulation of electromagnetic devices
described by nonlinear eddy current problems with harmonic excitations,
see, e.g.,
\cite{PhD:YamadaBessho:1988,
LRW:BachingerKaltenbacherReitzinger:2002a,
PhD:BachingerLangerSchoeberl:2005,
PhD:BachingerLangerSchoeberl:2006,
PhD:CopelandLanger:2010} and the references therein.
Later, the MhFEM has been applied to linear time-periodic parabolic boundary value 
and optimal control problems \cite{PhD:KollmannKolmbauer:2011,
LRW:KollmannKolmbauerLangerWolfmayrZulehner:2013, PhD:KrendlSimonciniZulehner:2012,
LRW:LangerWolfmayr:2013, LRW:Wolfmayr:2014}
and to linear time-periodic eddy current problems
and the corresponding optimal control problems
\cite{PhD:Kolmbauer:2012c, PhD:KolmbauerLanger:2012, PhD:KolmbauerLanger:2013}.
In the MhFEM setting, we are able to establish inf-sup and sup-sup conditions, which
provide existence and uniqueness of the solution to parabolic time-periodic problems.
For linear time-periodic parabolic problems, MhFEM
is a natural and very efficient numerical technology 
based on the decoupling of computations related to different modes.

This paper is aimed to make a step towards the creation of fully
reliable error estimation methods for distributed time-periodic parabolic 
optimal control problems.
We consider the multiharmonic finite element (MhFE) approximations of the 
reduced optimality system,
and derive
guaranteed and fully computable bounds for the discretization  errors.
For this purpose, we use the functional {\it a posteriori} error estimation
techniques earlier introduced by S.~Repin, see, e.g.,
the papers on parabolic problems
\cite{PhD:Repin:2002,
PhD:GaevskayaRepin:2005}
as well as on optimal control problems
\cite{PhD:GaevskayaHoppeRepin:2006,
PhD:GaevskayaHoppeRepin:2007},
the books
\cite{PhD:Repin:2008,
PhD:MaliNeittaanmaekiRepin:2014}
and the references therein.
In particular, our functional {\it a posteriori} error analysis
uses the techniques close to those suggested
in \cite{PhD:Repin:2002},
but the analysis contains essential changes due to the MhFEM setting.
In \cite{LRW:LangerRepinWolfmayr:2015}, the authors already derived 
functional {\it a posteriori} error estimates for MhFE approximations
to parabolic time-periodic boundary value problems.
Similar results are now obtained for the MhFE approximations
to the state and co-state, which are the unique solutions of the  
reduced optimality system. 
It is worth mentioning
that these {\it a posteriori} error estimates for the state and co-state 
immediately yield the corresponding {\it a posteriori} error estimates for the control. 
In addition to these results, we deduce
fully computable estimates of the cost functional.
In fact, we generate a new formulation of the optimal control
problem, in which (unlike the original statement) 
the state equations are accounted in terms of penalties. 
It is proved that the modified cost functional attains its infimum at
the optimal control and the respective state function.
Therefore, in principle, we can use it as
an object of direct minimization, which value on each step provides a guaranteed
upper bound of the cost functional.

The paper is organized as follows:
In Section~\ref{Sec2:ParTimePerOCP}, we discuss
the time-periodic parabolic optimal control problem and the corresponding
optimality system.
The multiharmonic finite element discretization of this space-time weak formulation
is considered in Section~\ref{Sec3:MhFEApprox}.
Section~\ref{Sec4:FunctionalAPostErrorEstimates:OptSys} is devoted to
the derivation of functional {\it a posteriori} error
estimates for the optimality system associated with the multiharmonic setting.
In Section~\ref{Sec5:APosterioriErrorEstimation:OCP:CostFuncs}, 
we present new results related to 
guaranteed and computable bounds of the cost functional.
In the final Section~\ref{Sec6:NumericalResults}, we discuss some implementation
issues and present the first numerical results.

\section{A Time-Periodic Parabolic Optimal Control Problem}
\label{Sec2:ParTimePerOCP}

Let $Q_T := \Omega \times (0,T)$ denote the space-time cylinder
and $\Sigma_T := \Gamma \times (0,T)$ its mantle boundary,
where the spatial domain $\Omega \subset \mathbb{R}^d$, $d=\{1,2,3\}$, is assumed to be
a bounded Lipschitz domain with boundary $\Gamma := \partial \Omega$,
and $(0,T)$ is a given time interval.
Let $\lambda>0$ be the regularization or cost parameter.
We consider the following parabolic time-periodic optimal control problem:
\begin{align}
\label{equation:minfunc:OCP}
 \min_{y,u} \mathcal{J}(y,u)
	:=  \frac{1}{2} \int_0^T \int_{\Omega} 
	\left(y(\boldsymbol{x}, t) - y_d(\boldsymbol{x},t) \right)^2 d\boldsymbol{x}\,dt
	+ \frac{\lambda}{2} \int_0^T \int_{\Omega} \left( u(\boldsymbol{x},t) \right)^2 d\boldsymbol{x}\,dt
\end{align}
subject to the parabolic time-periodic boundary value problem
\begin{align}
\left.
\label{equation:forwardpde:OCP}
\begin{aligned}
 \sigma(\boldsymbol{x}) \, \partial_t y(\boldsymbol{x},t)
 - \text{div} \, ( \nu(\boldsymbol{x}) \nabla y(\boldsymbol{x},t))
 &= u(\boldsymbol{x},t) \hspace{1cm} &(\boldsymbol{x},t) \in Q_T, \\
 y(\boldsymbol{x},t) &= 0 \hspace{1cm} &(\boldsymbol{x},t) \in \Sigma_T, \\
 y(\boldsymbol{x},0) &= y(\boldsymbol{x},T) \hspace{1cm} &\boldsymbol{x} \in \overline{\Omega}.
\end{aligned}
\quad
\right \rbrace
\end{align}
Here, $y$ and $u$ are the state and control functions, respectively.
The coefficients $\sigma(\cdot)$ and $\nu(\cdot)$
are uniformly bounded and satisfy the conditions
\begin{align}
 \label{assumptions:sigmaNu:sigmaStrictlyPositive}
 0 < \underline{\sigma} \leq \sigma(\boldsymbol{x}) \leq \overline{\sigma}, \qquad
 \text{and} \qquad
 0 < \underline{\nu} \leq \nu(\boldsymbol{x}) \leq \overline{\nu}, \qquad 
 \boldsymbol{x} \in \Omega,
\end{align}
where $\underline{\sigma}$, $\overline{\sigma}$, $\underline{\nu}$ and, $\overline{\nu}$ are constants.
As usual, the cost functional (\ref{equation:minfunc:OCP}) contains
a penalty term weighted with a positive factor $\lambda$. This term restricts
(in the integral sense) values of the control function $u$.

We can reformulate the problem
(\ref{equation:minfunc:OCP})-(\ref{equation:forwardpde:OCP})
in an equivalent form.
For this purpose, we introduce the Lagrangian
\begin{align}
 \label{equation:LagrangeFunctional}
 \mathcal{L}(y,u,p):= \mathcal{J}(y,u) - \int_0^T \int_{\Omega} \big(\sigma \partial_t y
    - \text{div} \, (\nu \nabla y) - u\big) p \, d\boldsymbol{x}\,dt,
\end{align}
which has a saddle point, see, e.g., \cite{LRW:HinzePinnauUlbrichUlbrich:2009, LRW:Troeltzsch:2010}
and the references therein.
The proper sets for which this saddle point problem is correctly stated are defined later.
Since the saddle point exists, 
the corresponding solutions satisfy the system of necessary conditions
\begin{align}
\label{equation:optsysParabolicTimePeriodicOCP}
\begin{aligned}
 \mathcal{L}_p (y,u,p) = 0, \qquad \qquad
 \mathcal{L}_y (y,u,p) = 0, \qquad \qquad
 \mathcal{L}_u (y,u,p) = 0.
\end{aligned}
\end{align}
Using the second condition, we can eliminate the control $u$ from the optimality system
(\ref{equation:optsysParabolicTimePeriodicOCP}), i.e.,
\begin{align}
 \label{equation:eliminateControl}
	u= - \lambda^{-1} p\; \;\mbox{ in}\;\; Q_T.
\end{align}
From (\ref{equation:eliminateControl}) it appears very natural to choose
$u$ and $p$ from the same space.
Moreover, we arrive at a reduced optimality system, written in its classical
formulation as
\begin{align}
\label{equation:KKTSysClassical}
\left.
\begin{aligned}
   \sigma(\boldsymbol{x}) \, \partial_t y(\boldsymbol{x},t) - 
   \text{div} \, (\nu(\boldsymbol{x}) \nabla y(\boldsymbol{x},t))
	&= - \lambda^{-1} p(\boldsymbol{x},t) 
	\quad &(\boldsymbol{x},t) \in Q_T, \\
   y(\boldsymbol{x},t) &= 0 
	\quad &(\boldsymbol{x},t) \in \Sigma_T, \\
   y(\boldsymbol{x},0) &= y(\boldsymbol{x},T) 
	\quad &\boldsymbol{x} \in \overline{\Omega},\\
   - \sigma(\boldsymbol{x}) \, \partial_t p(\boldsymbol{x},t) 
   - \text{div} \, (\nu(\boldsymbol{x}) \nabla p(\boldsymbol{x},t))
	&= y(\boldsymbol{x},t) - y_d(\boldsymbol{x},t) 
	\quad &(\boldsymbol{x},t) \in Q_T ,\\
  p(\boldsymbol{x},t) &= 0 
	\quad &(\boldsymbol{x},t) \in \Sigma_T ,\\
  p(\boldsymbol{x},T) &= p(\boldsymbol{x},0) 
	\quad &\boldsymbol{x} \in \overline{\Omega}.
\end{aligned}
\quad
\right \rbrace
\end{align}
In order to determine generalized solutions
of (\ref{equation:KKTSysClassical}),
we define the following
spaces (here and later on, the notation is similar to that was used in
\cite{LRW:Ladyzhenskaya:1973, LRW:LadyzhenskayaSolonnikovUralceva:1968}):
\begin{align*}
 H^{1,0}(Q_T) &= \{v \in L^2(Q_T) : \nabla v \in [L^2(Q_T)]^d \}, \\
 H^{0,1}(Q_T) &= \{v \in L^2(Q_T) : \partial_t v \in L^2(Q_T) \}, \\
 H^{1,1}(Q_T) &= \{v \in L^2(Q_T) : \nabla v \in [L^2(Q_T)]^d, \partial_t v \in L^2(Q_T) \},
\end{align*}
which are endowed with the norms
\begin{align*}
 \|v\|_{1,0} &:= \left(\int_{Q_T}
			\left(v(\boldsymbol{x},t)^2 + |\nabla v(\boldsymbol{x},t)|^2 \right)
			\, d\boldsymbol{x} \, dt \right)^{1/2}, \\
 \|v\|_{0,1} &:= \left(\int_{Q_T} \left(v(\boldsymbol{x},t)^2 + |\partial_t v(\boldsymbol{x},t)|^2 \right) \, d\boldsymbol{x} \, dt \right)^{1/2}, \\
 \|v\|_{1,1} &:= \left(\int_{Q_T} \left(v(\boldsymbol{x},t)^2 + |\nabla v(\boldsymbol{x},t)|^2
			+ |\partial_t v(\boldsymbol{x},t)|^2 \right) \, d\boldsymbol{x} \, dt \right)^{1/2},
\end{align*}
respectively. Here, 
$\nabla = \nabla_{\boldsymbol{x}}$ is the spatial gradient and $\partial_t$ denotes the generalized
derivative with respect to time. 
Also we define subspaces of the above introduced spaces by putting subindex zero if the functions
satisfy the homogeneous Dirichlet condition on $\Sigma_T$ and subindex \textit{per} if they
satisfy the periodicity condition $v(\boldsymbol{x},0) = v(\boldsymbol{x},T)$.
All inner products and norms in $L^2$ related to the whole space-time domain $Q_T$
are denoted by $\langle \cdot, \cdot \rangle$ and $\| \cdot \|$, respectively.
If they are associated with the spatial domain $\Omega$, then we write
$\langle \cdot, \cdot \rangle_{\Omega}$ and $\| \cdot \|_{\Omega}$. 
The symbols $\langle \cdot, \cdot \rangle_{1,\Omega}$ and $\| \cdot \|_{1,\Omega}$
denote the standard inner products and norms of the space $H^1(\Omega)$. 

The functions used in our analysis are presented in 
terms of Fourier series, e.g., for the function $v(\boldsymbol{x},t)$
is representation is
\begin{align}
\label{def:FourierAnsatz}
 v(\boldsymbol{x},t) = v_0^c(\boldsymbol{x}) +
      \sum_{k=1}^{\infty} \left(v_k^c(\boldsymbol{x}) \cos(k \omega t)
                          + v_k^s(\boldsymbol{x}) \sin(k \omega t)\right).
\end{align}
Here,
\begin{align*}
 \begin{aligned}
 v_0^c(\boldsymbol{x}) &= \frac{1}{T} \int_0^T v(\boldsymbol{x},t) \,dt, \\
 v_k^c(\boldsymbol{x}) = \frac{2}{T} \int_0^T v(\boldsymbol{x},t) \cos(k \omega t)\,&dt, 
 \hspace{0.3cm} \mbox{ and } \hspace{0.3cm}
 v_k^s(\boldsymbol{x}) = \frac{2}{T} \int_0^T v(\boldsymbol{x},t) \sin(k \omega t)\,dt 
 \end{aligned}
\end{align*}
are the Fourier coefficients, 
$T$ denotes the periodicity and $\omega = 2 \pi /T$ is the frequency. 
Since the problem has time-periodical conditions,
these representations of the exact solution and respective 
approximations are quite natural.
In what follows, we also use the spaces 
\begin{align*}
 H^{0,\frac{1}{2}}_{per}(Q_T) &= \{ v \in L^2(Q_T) : \big\| \partial^{1/2}_t v \big\|  < \infty \}, \\
 H^{1,\frac{1}{2}}_{per}(Q_T) &= \{ v \in H^{1,0}(Q_T) : \big\| \partial^{1/2}_t v \big\|  < \infty \}, \\
 H^{1,\frac{1}{2}}_{0,per}(Q_T) &= \{ v \in H^{1,\frac{1}{2}}_{per}(Q_T): v = 0 \mbox{ on } \Sigma_T \},
\end{align*}
where $\big\| \partial^{1/2}_t v \big\|$ is defined
in the Fourier space by the relation
\begin{align}
\label{definition:H01/2seminorm}
 \big\| \partial^{1/2}_t v \big\| ^2 :=
 |v|_{0,\frac{1}{2}}^2 := 
 \frac{T}{2} \sum_{k=1}^{\infty} k \omega \|\boldsymbol{v}_k\|_{\Omega}^2.
\end{align}
Here, $\boldsymbol{v}_k = (v_k^c,v_k^s)^T$ for all $k \in \mathbb{N}$,
see also \cite{LRW:LangerWolfmayr:2013}. 
These spaces can be considered as Hilbert spaces if we introduce the following
(equivalent) products:
\begin{align*}
\begin{aligned}
 \langle \partial^{1/2}_t y, \partial^{1/2}_t v \rangle  :=
 \frac{T}{2} \sum_{k=1}^{\infty} k \omega \langle\boldsymbol{y}_k,\boldsymbol{v}_k\rangle_{\Omega},
 \qquad \quad
 \langle\sigma
 \partial^{1/2}_t y, \partial^{1/2}_t v\rangle  :=
 \frac{T}{2} \sum_{k=1}^{\infty} k \omega \langle\sigma
 \boldsymbol{y}_k,\boldsymbol{v}_k\rangle_{\Omega}.
\end{aligned}
\end{align*}
The above introduced spaces allow us to operate with a "symmetrized" formulation of the problem
(\ref{equation:KKTSysClassical}) presented by (\ref{problem:KKTSysSTVFAPost}).
The seminorm and the norm of the space $H^{1,\frac{1}{2}}_{per}(Q_T)$ are defined 
by the relations
\begin{align*}
 |v|_{1,\frac{1}{2}}^2 
 &:= \|\nabla v\| ^2 + \|\partial_t^{1/2} v\| ^2 
 = T \, \|\nabla y_0^c\|_{\Omega}^2
 + \frac{T}{2} \sum_{k=1}^{\infty} \left(k\omega \|\boldsymbol{v}_k\|_{\Omega}^2
 + \|\nabla \boldsymbol{v}_k\|_{\Omega}^2\right)
\qquad \mbox{ and } \\
 \|v\|_{1,\frac{1}{2}}^2
 &:= \|v\| ^2 + |v|_{1,\frac{1}{2}}^2 
 = T \, (\|y_0^c\|_{\Omega}^2 + \|\nabla y_0^c\|_{\Omega}^2)
 + \frac{T}{2} \sum_{k=1}^{\infty} \left((1+k\omega) \|\boldsymbol{v}_k\|_{\Omega}^2
 + \|\nabla \boldsymbol{v}_k\|_{\Omega}^2\right),
\end{align*}
respectively.
Using Fourier type series, 
it is easy to define the function
"orthogonal" to $v$: 
\begin{align*}
 \begin{aligned}
 v^{\perp}(\boldsymbol{x},t) &:= 
 \sum_{k=1}^{\infty} 
 \left(- v_k^c(\boldsymbol{x}) \sin(k \omega t)
 + v_k^s(\boldsymbol{x}) \cos(k \omega t)\right) \\
 &= \sum_{k=1}^{\infty} 
 \underbrace{(v_k^s(\boldsymbol{x}),-v_k^c(\boldsymbol{x}))}_{=:(-\boldsymbol{v}_k^{\perp})^T} \cdot
                        \left( \begin{array}{l}
                          \cos(k \omega t) \\
                          \sin(k \omega t)
                        \end{array} \right).
 \end{aligned}
\end{align*}
Obviously, $\|\boldsymbol{u}_k^\perp\|_{\Omega} = \|\boldsymbol{u}_k\|_{\Omega}$ and
we find that 
\begin{align*}
 \big\|\partial^{1/2}_t v^\perp\big\| ^2
 = \frac{T}{2} \sum_{k=1}^\infty k \omega \|\boldsymbol{v}_k^\perp\|_{\Omega}^2
 = \frac{T}{2} \sum_{k=1}^\infty k \omega \|\boldsymbol{v}_k\|_{\Omega}^2
 = \big\|\partial^{1/2}_t v\big\| ^2 \qquad \forall \, v \in H^{0,\frac{1}{2}}_{per}(Q_T).
\end{align*}
Henceforth, we use the following subsidiary result (which 
proof can be found in \cite{LRW:LangerRepinWolfmayr:2015}):
\begin{lemma}
 \label{lemma:H11/2IdentiesAndOrthogonalities}
The identities
 \begin{align}
 \label{equation:H11/2identities}
 \begin{aligned}
  \langle \sigma \partial_t^{1/2} y,\partial_t^{1/2} v  \rangle  =
  \langle \sigma \partial_t y,v^{\perp}  \rangle  \quad \mbox{ and } \quad 
  \langle \sigma \partial_t^{1/2} y,\partial_t^{1/2} v^{\perp}  \rangle =
  \langle \sigma \partial_t y,v  \rangle
 \end{aligned}
 \end{align}
 are valid 
 for all $y \in H^{0,1}_{per}(Q_T)$ and $v \in H^{0,\frac{1}{2}}_{per}(Q_T)$.
\end{lemma}

Also, we recall the orthogonality relations (see \cite{LRW:LangerWolfmayr:2013, LRW:Wolfmayr:2014})
\begin{align}
\left.
\label{equation:orthorelation}
\begin{aligned}
\langle \sigma \partial_t y,y \rangle  = 0 \quad &\mbox{ and } \quad
 \langle \sigma y^{\perp},y \rangle  = 0 \qquad &&\forall \, y \in H^{0,1}_{per}(Q_T), \\
 \langle \sigma \partial^{1/2}_t y,\partial^{1/2}_t y^{\perp} \rangle  = 0 \quad &\mbox{ and } \quad
\langle \nu \nabla y, \nabla y^{\perp} \rangle  = 0 \qquad &&\forall \, y \in H^{1,\frac{1}{2}}_{per}(Q_T),
\end{aligned}
\quad
\right \rbrace
\end{align}
and the identity 
\begin{align}
\label{def:identityH01/2per}
 \int_0^T \xi \, \partial_t^{1/2} v^\perp \, dt = - \int_0^T \partial_t^{1/2} \xi^\perp \, v \, dt
 \qquad \forall \, \xi, v \in H^{0,\frac{1}{2}}_{per}(Q_T),
\end{align}
where
\begin{align}
 \langle \xi,\partial_t^{1/2} v \rangle
 := \frac{T}{2} \sum_{k=1}^{\infty} (k \omega)^{1/2} \langle \boldsymbol{\xi}_k,\boldsymbol{v}_k \rangle_{\Omega}.
\end{align}

We note that for functions presented in terms of Fourier series the standard
Friedrichs inequality holds in $Q_T$. Indeed,
\begin{align}
\label{inequality:Friedrichs:FourierSpace}
 \begin{aligned}
 \|\nabla u\| ^2 &= \int_{Q_T} |\nabla u|^2 \, d\boldsymbol{x}\,dt
 = T \, \|\nabla u_0^c\|_{\Omega}^2 + \frac{T}{2} \sum_{k=1}^\infty \|\nabla \boldsymbol{u}_k\|_{\Omega}^2 \\
 &\geq \frac{1}{C_F^2} \left(T \, \|u_0^c\|_{\Omega}^2 + \frac{T}{2} \sum_{k=1}^\infty \|\boldsymbol{u}_k\|_{\Omega}^2 \right) 
 = \frac{1}{C_F^2} \|u\| ^2.
 \end{aligned}
\end{align}

In order to derive the weak formulations of the equations  
(\ref{equation:KKTSysClassical}), we multiply them
by the test functions $z, q \in H^{1,\frac{1}{2}}_{0,per}(Q_T)$,
integrate over $Q_T$, and apply
integration by parts with respect to the spatial variables 
and time. 
We arrive at
the following ``symmetric'' space-time weak 
formulations of the reduced
optimality system (\ref{equation:KKTSysClassical}):
Given the desired state $y_d \in L^2(Q_T)$, find $y$ and $p$ from 
$H^{1,\frac{1}{2}}_{0,per}(Q_T)$ such that
\begin{align}
 \label{problem:KKTSysSTVFAPost}
 \left.
 \begin{aligned}
  &\int_{Q_T} \Big( y\,z - \nu(\boldsymbol{x}) \nabla p \cdot \nabla z
      + \sigma(\boldsymbol{x}) \partial_t^{1/2} p \, \partial_t^{1/2} z^{\perp} \Big)\,d\boldsymbol{x}\,dt
      = \int_{Q_T} y_d\,z\,d\boldsymbol{x}\,dt, \\
  &\int_{Q_T} \Big( \nu(\boldsymbol{x}) \nabla y \cdot \nabla q
      + \sigma(\boldsymbol{x}) \partial_t^{1/2} y \, \partial_t^{1/2} q^{\perp}
      + \lambda^{-1} p\,q \Big)\,d\boldsymbol{x}\,dt = 0,
 \end{aligned}
\quad
\right \rbrace
\end{align}
for all test functions $z, q \in H^{1,\frac{1}{2}}_{0,per}(Q_T)$.
We can represent (\ref{problem:KKTSysSTVFAPost}) in a somewhat different form.
For this purpose, it is convenient to introduce the bilinear form 
\begin{align}
\label{definition:KKTSysSTVF}
\begin{aligned}
  \mathcal{B}((y,p),(z,q)) = \int_{Q_T} &\Big( y\,z - \nu(\boldsymbol{x}) \nabla p \cdot \nabla z
      + \sigma(\boldsymbol{x}) \partial_t^{1/2} p \, \partial_t^{1/2} z^{\perp} \\
      &+ \nu(\boldsymbol{x}) \nabla y \cdot \nabla q
      + \sigma(\boldsymbol{x}) \partial_t^{1/2} y \, \partial_t^{1/2} q^{\perp}
      + \lambda^{-1} p\,q  \Big)\,d\boldsymbol{x}\,dt.
\end{aligned}
\end{align}
Then (\ref{problem:KKTSysSTVFAPost}) reads as
\begin{align*}
 \mathcal{B}((y,p),(z,q)) = \langle (y_d,0),(z,q)\rangle \qquad \forall \, (z, q) \in (H^{1,\frac{1}{2}}_{0,per}(Q_T))^2.
\end{align*}

\section{Multiharmonic Finite Element Approximation}
\label{Sec3:MhFEApprox}

In order to solve the 
optimal control problem (\ref{equation:minfunc:OCP})-(\ref{equation:forwardpde:OCP}) approximately, 
we discretize the optimality system (\ref{problem:KKTSysSTVFAPost}) by
the MhFEM (see \cite{LRW:LangerWolfmayr:2013}).
Using the Fourier series ansatz (\ref{def:FourierAnsatz}) in
(\ref{problem:KKTSysSTVFAPost}) and exploiting 
the orthogonality of $\cos(k \omega t)$ and $\sin(k \omega t)$,
we arrive at the following problem:
Find $\boldsymbol{y}_k, \boldsymbol{p}_k \in \mathbb{V} := V \times V = (H^1_0(\Omega))^2$ such that
\begin{align}
 \label{equation:MultiAnsVFBlockk}
 \left.
\begin{aligned}
 &\int_{\Omega} \big( \boldsymbol{y}_k \cdot \boldsymbol{z}_k
  - \nu(\boldsymbol{x}) \nabla \boldsymbol{p}_k \cdot \nabla \boldsymbol{z}_k
  + k \omega \, \sigma(\boldsymbol{x}) \boldsymbol{p}_k \cdot \boldsymbol{z}_k^{\perp} \big)\,d\boldsymbol{x}
  = \int_{\Omega} \boldsymbol{y_d}_k \cdot \boldsymbol{z}_k\,d\boldsymbol{x}, \\
 &\int_{\Omega} \big( \nu(\boldsymbol{x}) \nabla \boldsymbol{y}_k \cdot \nabla \boldsymbol{q}_k
  + k \omega \, \sigma(\boldsymbol{x}) \boldsymbol{y}_k \cdot \boldsymbol{q}_k^{\perp}
  + \lambda^{-1} \boldsymbol{p}_k \cdot \boldsymbol{q}_k \big)\,d\boldsymbol{x} = 0,
\end{aligned}
\quad
\right \rbrace
\end{align}
for all test functions $\boldsymbol{z}_k, \boldsymbol{q}_k \in \mathbb{V}$. 
The system (\ref{equation:MultiAnsVFBlockk}) must be solved for every
mode $k \in \mathbb{N}$.
For $k = 0$, we obtain a reduced problem:
Find $y_0^c, p_0^c \in V$ such that
\begin{align}
 \label{equation:MultiAnsVFBlock0}
 \left.
\begin{aligned}
 &\int_{\Omega} \big( y_0^c \cdot z_0^c - \nu(\boldsymbol{x}) \nabla p_0^c \cdot \nabla z_0^c \big)\,d\boldsymbol{x}
 = \int_{\Omega} {y_d^c}_0 \cdot z_0^c\,d\boldsymbol{x}, \\
 &\int_{\Omega} \big( \nu(\boldsymbol{x}) \nabla y_0^c \cdot \nabla q_0^c
 + \lambda^{-1} p_0^c \cdot q_0^c \big)\,d\boldsymbol{x} = 0,
\end{aligned}
\quad
\right \rbrace
\end{align}
for all test functions $z_0^c, q_0^c \in V$. 
The problems (\ref{equation:MultiAnsVFBlockk}) and (\ref{equation:MultiAnsVFBlock0})
have unique solutions.
In order to solve these problems numerically, the Fourier series are truncated at a finite index $N$ and
the unknown Fourier coefficients
$ \boldsymbol{y}_k = (y_k^c, y_k^s)^T, \, \boldsymbol{p}_k = (p_k^c, p_k^s)^T \in \mathbb{V} $
are approximated by finite element (FE) functions
\begin{align*}
 \boldsymbol{y}_{kh} = (y_{kh}^c, y_{kh}^s)^T, \, \boldsymbol{p}_{kh} = (p_{kh}^c, p_{kh}^s)^T
\in \mathbb{V}_h = V_h \times V_h \subset \mathbb{V},
\end{align*}
where $ V_h = \mbox{span} \{\varphi_1, \dots, \varphi_n\} $
with $\{\varphi_i(\boldsymbol{x}): i=1,2,\dots,n_h \}$ 
is a conforming FE space. 
We denote by $h$ the usual discretization parameter such that $n = n_h = \mbox{dim} V_h = O(h^{-d})$.
In this work, we  
use continuous, piecewise linear finite elements on 
a regular triangulation $\mathcal{T}_h$ to construct 
$V_h$ and its basis (see, e.g.,
\cite{LRW:Braess:2005, LRW:Ciarlet:1978, LRW:JungLanger:2013, LRW:Steinbach:2008}).
This leads to the following saddle point system
for every single mode $k=1,2,\dots,N$:
\begin{align}
 \label{equation:MultiFESysBlockk}
 \left( \begin{array}{cccc}
     M_h  &  0 & -K_{h,\nu} & k \omega M_{h,\sigma} \\
     0  &  M_h & -k \omega M_{h,\sigma} & -K_{h,\nu} \\
     -K_{h,\nu}  &  -k \omega M_{h,\sigma} & -\lambda^{-1} M_h & 0 \\
     k \omega M_{h,\sigma}  &  -K_{h,\nu} & 0 & -\lambda^{-1} M_h \end{array} \right) \left( \begin{array}{c}
     \underline{y}_k^c \\
     \underline{y}_k^s \\
     \underline{p}_k^c \\
     \underline{p}_k^s \end{array} \right) = \left( \begin{array}{c}
     {\underline{y}_d^c}_k \\
     {\underline{y}_d^s}_k \\
     0 \\
     0 \end{array} \right),
\end{align}
which has to be solved with respect to the nodal parameter vectors
\begin{align*}
\underline{y}_k^c = ( y_{k,i}^c)_{i=1,\dots,n}, \,
\underline{y}_k^s = ( y_{k,i}^s)_{i=1,\dots,n}, \, 
\underline{p}_k^c = ( p_{k,i}^c)_{i=1,\dots,n}, \,
\underline{p}_k^s = ( p_{k,i}^s)_{i=1,\dots,n} \in \mathbb{R}^n
\end{align*}
of the FE approximations 
$y _{kh}^c(\boldsymbol{x}) = \sum_{i=1}^n y_{k,i}^c \varphi_i(\boldsymbol{x})$ 
and
$y _{kh}^s(\boldsymbol{x}) = \sum_{i=1}^n y_{k,i}^s \varphi_i(\boldsymbol{x})$.  
Similarly, 
$p _{kh}^c(\boldsymbol{x}) = \sum_{i=1}^n p_{k,i}^c \varphi_i(\boldsymbol{x})$ 
and
$p _{kh}^s(\boldsymbol{x}) = \sum_{i=1}^n p_{k,i}^s \varphi_i(\boldsymbol{x})$.
The matrices $M_h$, $M_{h,\sigma}$, and $K_{h,\nu}$ denote the mass matrix,
the weighted mass matrix and the stiffness matrix, respectively. Their entries are
defined by the integrals
\begin{align*}
\begin{aligned}
 M_h^{ij} = \int_{\Omega} \varphi_i \varphi_j \,d\boldsymbol{x}, \hspace{0.5cm}
 M_{h,\sigma}^{ij} &= \int_{\Omega} \sigma \, \varphi_i \varphi_j \,d\boldsymbol{x}, \hspace{0.5cm}
 K_{h,\nu}^{ij} &= \int_{\Omega} \nu \, \nabla \varphi_i \cdot \nabla \varphi_j \,d\boldsymbol{x} 
\end{aligned}
\end{align*}
and the right hand side vectors have the form
\begin{align*}
\begin{aligned}
 {\underline{y}_d^c}_k = \Big\lbrack \int_{\Omega} {y_d^c}_k \varphi_j \,d\boldsymbol{x}
 \Big\rbrack_{j=1,\dots,n} \quad \mbox{and} \quad
 {\underline{y}_d^s}_k = \Big\lbrack \int_{\Omega} {y_d^s}_k \varphi_j \,d\boldsymbol{x}
 \Big\rbrack_{j=1,\dots,n}.
\end{aligned}
\end{align*}
For $k=0$, the problem (\ref{equation:MultiAnsVFBlock0}) generates a reduced 
system of linear equations, i.e., 
\begin{align}
 \label{equation:MultiFESysBlock0}
 \left( \begin{array}{cc}
     M_h  &  -K_{h,\nu} \\
     -K_{h,\nu}  &  - \lambda^{-1} M_h \end{array} \right) \left( \begin{array}{c}
     \underline{y}_0^c \\
     \underline{p}_0^c \end{array} \right) = \left( \begin{array}{c}
     {\underline{y}_d^c}_0 \\
     0 \end{array} \right).
\end{align}
Fast and robust solvers for the systems
(\ref{equation:MultiFESysBlockk}) and (\ref{equation:MultiFESysBlock0})
can be found in 
\cite{LRW:KollmannKolmbauerLangerWolfmayrZulehner:2013, LRW:KrausWolfmayr:2013,
LRW:LangerWolfmayr:2013, LRW:Wolfmayr:2014},
which we use in order to obtain the MhFE approximations
\begin{align}
\left.
\label{definition:MultiharmonicFEAnsatzStateANDAdjointState}
\begin{aligned}
 y_{N h}(\boldsymbol{x},t) &= y_{0h}^c(\boldsymbol{x}) + \sum_{k=1}^N \left(y_{kh}^c(\boldsymbol{x}) \cos(k \omega t)
      + y_{kh}^s(\boldsymbol{x}) \sin(k \omega t)\right), \\
 p_{N h}(\boldsymbol{x},t) &= p_{0h}^c(\boldsymbol{x}) + \sum_{k=1}^N \left(p_{kh}^c(\boldsymbol{x}) \cos(k \omega t)
      + p_{kh}^s(\boldsymbol{x}) \sin(k \omega t)\right).
\end{aligned}
\quad
\right \rbrace
\end{align}

\section{Functional A Posteriori Error Estimates for the Optimality System}
\label{Sec4:FunctionalAPostErrorEstimates:OptSys}

Now we are concerned with {\it a posteriori} estimates of the difference
between the exact solution $(y,p)$
and the respective finite element solution $(y_{N h},p_{N h})$.
First, we present the inf-sup and sup-sup conditions
for the bilinear form (\ref{definition:KKTSysSTVF}).
\begin{lemma}
\label{lemma:KKTSysSTVFinfsupsupsup}
  For all $y,p \in H^{1,\frac{1}{2}}_{0,per}(Q_T)$,
 the space-time bilinear form 
 $\mathcal{B}(\cdot,\cdot)$ defined by
 (\ref{definition:KKTSysSTVF})
 meets the following inf-sup and sup-sup conditions:
 \begin{align}
 \label{inequality:KKTSysSTVFinfsupsupsup}
  \mu_{1} \|(y,p)\|_{1,\frac{1}{2}}
  \leq
  \sup_{0 \not= (z,q) \in (H^{1,\frac{1}{2}}_{0,per}(Q_T))^2}
    \frac{\mathcal{B}((y,p),(z,q))}{\|(z,q)\|_{1,\frac{1}{2}}}
    \leq
  \mu_{2} \|(y,p)\|_{1,\frac{1}{2}}, 
 \end{align}
where 
 $\mu_1 = \frac{\min\{\frac{1}{\sqrt{\lambda}},\underline{\nu},\underline{\sigma}\}
 \min\{\sqrt{\lambda},\frac{1}{\sqrt{\lambda}}\}}{\sqrt{1+2\max\{\lambda,\frac{1}{\lambda}\}}}$
 and
 $\mu_{2} = \max\{1,\frac{1}{\lambda},\overline{\nu},\overline{\sigma}\}$
 are positive constants.
\begin{proof}
Using the triangle and Cauchy-Schwarz inequalities, we obtain 
\begin{align*}
  \big|\mathcal{B}((y,p),(z,q))\big|
  = &\, \Big| \int_0^T \int_\Omega \Big( y\,z - \nu(\boldsymbol{x}) \nabla p \cdot \nabla z
      + \sigma(\boldsymbol{x}) \partial_t^{1/2} p \, \partial_t^{1/2} z^{\perp} \\
      &\qquad \qquad \, + \nu(\boldsymbol{x}) \nabla y \cdot \nabla q
      + \sigma(\boldsymbol{x}) \partial_t^{1/2} y \, \partial_t^{1/2} q^{\perp}
      + \lambda^{-1} p\,q  \Big)\,d\boldsymbol{x}\,dt \Big| \\
  \leq &\, \|y\| \|z\|
      + \overline{\nu} \, \|\nabla p\|\|\nabla z\|
      + \overline{\sigma} \, \big\|\partial^{1/2}_t p\big\| \big\|\partial^{1/2}_t z\big\| \\
      &+ \overline{\nu} \, \|\nabla y\| \|\nabla q\|
      + \overline{\sigma} \, \big\|\partial^{1/2}_t y\big\| \big\|\partial^{1/2}_t q\big\|
      + \lambda^{-1} \, \|p\| \|q\| \\
   \leq 
  &\, \mu_2 \|(y,p)\|_{1,\frac{1}{2}} \|(z,q)\|_{1,\frac{1}{2}},
\end{align*}
where
$\mu_2 := \max\{1,\frac{1}{\lambda},\overline{\nu},\overline{\sigma}\}$.
Thus, the right hand-side inequality in (\ref{inequality:KKTSysSTVFinfsupsupsup})
is proved.
In order to prove the left hand-side inequality, we select the test functions
\begin{align*}
 (z,q) = (y - \frac{1}{\sqrt{\lambda}} p - \frac{1}{\sqrt{\lambda}} p^\perp,
          p + \sqrt{\lambda} y - \sqrt{\lambda} y^\perp).
\end{align*}
Using the $\sigma$- and $\nu$-weighted orthogonality relations (\ref{equation:orthorelation}),
we obtain the relations
\begin{align*}
 \mathcal{B}((y,p),(y,p)) 
  =&\, \|y\|^2 + \lambda^{-1} \|p\|^2, \\
 \mathcal{B}((y,p),(-\frac{1}{\sqrt{\lambda}} p, \sqrt{\lambda} y)) 
  = &\, \frac{1}{\sqrt{\lambda}} \langle \nu \nabla p,\nabla p \rangle
  + \sqrt{\lambda} \langle \nu \nabla y,\nabla y \rangle, \\
 \mathcal{B}((y,p),(-\frac{1}{\sqrt{\lambda}} p^\perp,-\sqrt{\lambda} y^\perp)) 
  = &\, \frac{1}{\sqrt{\lambda}} \langle \sigma \partial_t^{1/2} p, \partial_t^{1/2} p \rangle
    + \sqrt{\lambda} \langle \sigma \partial_t^{1/2} y, \partial_t^{1/2} y \rangle,
\end{align*}
which lead to the estimate
\begin{align*}
 \mathcal{B}((y,p)&,(y - \frac{1}{\sqrt{\lambda}} p - \frac{1}{\sqrt{\lambda}} p^\perp,
          p + \sqrt{\lambda} y - \sqrt{\lambda} y^\perp)) \\
 &\geq 
 \min\{\frac{1}{\sqrt{\lambda}},\underline{\nu},\underline{\sigma}\}
 \min\{\sqrt{\lambda},\frac{1}{\sqrt{\lambda}}\}
\|(y,p)\|_{1,\frac{1}{2}}^2.
\end{align*}
Since
\begin{align*}
 \|(z,q)\|_{1,\frac{1}{2}}^2 
 \leq \left(1+2\max\{\lambda,\frac{1}{\lambda}\}\right) \|(y,p)\|_{1,\frac{1}{2}}^2,
\end{align*}
we arrive at the estimate 
\begin{align*}
  \sup_{0 \not= (z,q) \in (H^{1,\frac{1}{2}}_{0,per}(Q_T))^2}
    \frac{\mathcal{B}((y,p),(z,q))}{\|(z,q)\|_{1,\frac{1}{2}}}
   &\geq \frac{\min\{\frac{1}{\sqrt{\lambda}},\underline{\nu},\underline{\sigma}\}
 \min\{\sqrt{\lambda},\frac{1}{\sqrt{\lambda}}\} \|(y,p)\|_{1,\frac{1}{2}}^2}{
 \sqrt{1+2\max\{\lambda,\frac{1}{\lambda}\}} \|(y,p)\|_{1,\frac{1}{2}}} \\
   &= \mu_1 \, \|(y,p)\|_{1,\frac{1}{2}},
\end{align*}
where $\mu_1 = \frac{\min\{\frac{1}{\sqrt{\lambda}},\underline{\nu},\underline{\sigma}\}
 \min\{\sqrt{\lambda},\frac{1}{\sqrt{\lambda}}\}}{\sqrt{1+2\max\{\lambda,\frac{1}{\lambda}\}}}$.
\end{proof}
\end{lemma}

In view of the Friedrichs inequality,
the norms $|\cdot|_{1,\frac{1}{2}}$ and 
$\|\cdot\|_{1,\frac{1}{2}}$ are equivalent. 
Therefore, Lemma \ref{lemma:KKTSysSTVFinfsupsupsup}  
implies the following result:
\begin{lemma}
\label{lemma:KKTSysSTVFinfsupsupsup:Seminorm}
For all $y,p \in H^{1,\frac{1}{2}}_{0,per}(Q_T)$, the bilinear form 
$\mathcal{B}(\cdot,\cdot)$ satisfies the inf-sup and sup-sup conditions 
 \begin{align}
 \label{inequality:KKTSysSTVFinfsupsupsup:Seminorm}
  \tilde \mu_1 |(y,p)|_{1,\frac{1}{2}}
  \leq
  \sup_{0 \not= (z,q) \in (H^{1,\frac{1}{2}}_{0,per}(Q_T))^2}
    \frac{\mathcal{B}((y,p),(z,q))}{|(z,q)|_{1,\frac{1}{2}}}
    \leq
  \tilde \mu_2 |(y,p)|_{1,\frac{1}{2}},
 \end{align}
  where 
  $\tilde \mu_1 = \frac{\min\{\underline{\nu},\underline{\sigma}\} \min\{\lambda,\frac{1}{\lambda}\}}{\sqrt{2}} > 0$
  and 
 $\tilde \mu_2 = \max\{1,\frac{1}{\lambda},\overline{\nu},\overline{\sigma}\} \max\{1,C_F^2+1\} > 0$, and
 $C_F$ is the Friedrichs constant.
\begin{proof}
The right hand-side inequality in (\ref{inequality:KKTSysSTVFinfsupsupsup:Seminorm})
results from the triangle and Cauchy-Schwarz inequalities and, the 
Friedrichs inequality (\ref{inequality:Friedrichs:FourierSpace}). Indeed,
\begin{align*}
  \big|\mathcal{B}((y,p),(z,q))\big|
  \leq \max\{1,\frac{1}{\lambda},\overline{\nu},&\,\overline{\sigma}\}
      \Big((C_F^2+1) \|\nabla y\|^2 + \big\|\partial^{1/2}_t y\big\|^2 \\
      &\qquad +  (C_F^2+1) \|\nabla p\|^2 + \big\|\partial^{1/2}_t p\big\|^2\Big)^{1/2} \\  
      &\times \Big((C_F^2+1)  \|\nabla z\|^2
      + \big\|\partial^{1/2}_t z\big\|^2 \\
      &\qquad + (C_F^2+1) \|\nabla q\|^2 + \big\|\partial^{1/2}_t q\big\|^2 \Big)^{1/2} \\ 
  \leq 
  \tilde \mu_2 |(y,p)|_{1,\frac{1}{2}} |&(z,q)|_{1,\frac{1}{2}},
\end{align*}
where $\tilde \mu_2 = \max\{1,\frac{1}{\lambda},\overline{\nu},\overline{\sigma}\} \max\{1,C_F^2+1\}$.

The left hand-side inequality in (\ref{inequality:KKTSysSTVFinfsupsupsup:Seminorm}) is 
proved quite similarly to the previous case. We select 
the test functions
\begin{align*}
 (z,q) = (- \frac{1}{\sqrt{\lambda}} p - \frac{1}{\sqrt{\lambda}} p^\perp,
          \sqrt{\lambda} y - \sqrt{\lambda} y^\perp),
\end{align*}
and use the $\sigma$- and $\nu$-weighted orthogonality relations
(\ref{equation:orthorelation}).
Then, we find that
\begin{align*}
 \mathcal{B}((y,p),(- \frac{1}{\sqrt{\lambda}} p - \frac{1}{\sqrt{\lambda}} p^\perp,
          \sqrt{\lambda} y - \sqrt{\lambda} y^\perp))
 &\geq \min\{\underline{\nu},\underline{\sigma}\} \min\{\sqrt{\lambda},\frac{1}{\sqrt{\lambda}}\}
  |(y,p)|_{1,\frac{1}{2}}^2.
\end{align*}
In view of the estimate
\begin{align*}
 |(z,q)|_{1,\frac{1}{2}}^2 
 \leq 2 \max\{\lambda,\frac{1}{\lambda}\} |(y,p)|_{1,\frac{1}{2}}^2,
\end{align*}
we obtain
\begin{align*}
  \sup_{0 \not= (z,q) \in (H^{1,\frac{1}{2}}_{0,per}(Q_T))^2}
    \frac{\mathcal{B}((y,p),(z,q))}{\|(z,q)\|_{1,\frac{1}{2}}}
   &\geq \frac{\min\{\underline{\nu},\underline{\sigma}\} \min\{\sqrt{\lambda},\frac{1}{\sqrt{\lambda}}\}
  |(y,p)|_{1,\frac{1}{2}}^2}{
 \sqrt{2} \max\{\sqrt{\lambda},\frac{1}{\sqrt{\lambda}}\} |(y,p)|_{1,\frac{1}{2}}} \\
   &= \tilde \mu_1 \, |(y,p)|_{1,\frac{1}{2}},
  \end{align*}
  where
$\tilde \mu_1 = \frac{\min\{\underline{\nu},\underline{\sigma}\} \min\{\lambda,\frac{1}{\lambda}\}}{\sqrt{2}}$.
This justifies the left hand-side inequality in (\ref{inequality:KKTSysSTVFinfsupsupsup:Seminorm}).
\end{proof}
\end{lemma}

Let $(\eta,\zeta)$ be an approximation of $(y,p)$, which is 
a bit more regular with respect to the time variable
than the exact solution $(y,p)$. 
Namely, we assume that $\eta, \zeta \in H^{1,1}_{0,per}(Q_T)$ 
(this assumption is
of course true for the MhFE approximations
$y_{N h}$ and $p_{N h}$ defined in (\ref{definition:MultiharmonicFEAnsatzStateANDAdjointState})).
Our goal is to deduce a computable upper bound of the error
$e := (y,p) - (\eta,\zeta)$
in $H^{1,\frac{1}{2}}_{0,per}(Q_T) \times H^{1,\frac{1}{2}}_{0,per}(Q_T)$. 
From (\ref{problem:KKTSysSTVFAPost}), it follows that the integral identity
\begin{align}
\left.
 \label{problem:KKTSysSTVFAPost:Error}
  \begin{aligned}
 \int_{Q_T} &\Big( (y-\eta)\,z - \nu(\boldsymbol{x}) \nabla (p-\zeta) \cdot \nabla z
      + \sigma(\boldsymbol{x}) \partial_t^{1/2} (p-\zeta) \, \partial_t^{1/2} z^{\perp} \\
      &+ \nu(\boldsymbol{x}) \nabla (y-\eta) \cdot \nabla q
      + \sigma(\boldsymbol{x}) \partial_t^{1/2} (y-\eta) \, \partial_t^{1/2} q^{\perp}
      + \lambda^{-1} (p-\zeta)\,q \Big)\,d\boldsymbol{x}\,dt \\
 = &\int_{Q_T} \Big(y_d\,z
      - \eta\,z + \nu(\boldsymbol{x}) \nabla \zeta \cdot \nabla z
      - \sigma(\boldsymbol{x}) \partial_t^{1/2} \zeta \, \partial_t^{1/2} z^{\perp} \\
      &- \nu(\boldsymbol{x}) \nabla \eta \cdot \nabla q
      - \sigma(\boldsymbol{x}) \partial_t^{1/2} \eta \, \partial_t^{1/2} q^{\perp}
      - \lambda^{-1} \zeta \,q \Big) \, d\boldsymbol{x}\,dt
 \end{aligned}
\quad
\right \rbrace
\end{align}
holds for all $z, q \in H^{1,\frac{1}{2}}_{0,per}(Q_T)$.
The linear functional
\begin{align*}
 \begin{aligned}
 \mathcal{F}_{(\eta,\zeta)}(z,q)
 = \int_{Q_T} &\Big(y_d\,z
      - \eta\,z + \nu(\boldsymbol{x}) \nabla \zeta \cdot \nabla z
      - \sigma(\boldsymbol{x}) \partial_t^{1/2} \zeta \, \partial_t^{1/2} z^{\perp} \\
      &- \nu(\boldsymbol{x}) \nabla \eta \cdot \nabla q
      - \sigma(\boldsymbol{x}) \partial_t^{1/2} \eta \, \partial_t^{1/2} q^{\perp}
      - \lambda^{-1} \zeta \,q \Big) \, d\boldsymbol{x}\,dt 
 \end{aligned}
\end{align*}
is defined on $(z,q) \in H^{1,\frac{1}{2}}_{0,per}(Q_T) \times H^{1,\frac{1}{2}}_{0,per}(Q_T)$.
It can be viewed as a quantity measuring the accuracy of 
(\ref{problem:KKTSysSTVFAPost:Error}) for any pair of test functions
$(z,q)$. Therefore, getting an upper bound of the error is reduced to the estimation of 
\begin{align}
 \label{inequality:supRHS:OCP}
 \begin{aligned}
 \sup_{0 \not= (z,q) \in (H^{1,\frac{1}{2}}_{0,per}(Q_T))^2}
    \frac{\mathcal{F}_{(\eta,\zeta)}(z,q)}{\|(z,q)\|_{1,\frac{1}{2}}} 
 \qquad \text{ or } \quad
 \sup_{0 \not= (z,q) \in (H^{1,\frac{1}{2}}_{0,per}(Q_T))^2}
    \frac{\mathcal{F}_{(\eta,\zeta)}(z,q)}{|(z,q)|_{1,\frac{1}{2}}}. 
 \end{aligned}
\end{align}
We reconstruct $\mathcal{F}_{(\eta,\zeta)}$ using 
the identity
\begin{align}
\label{equation:identityEtaH11/2}
 \langle \sigma \partial_t^{1/2} \eta,\partial_t^{1/2} z^{\perp} \rangle 
 = \langle \sigma \partial_t \eta,z \rangle 
 \qquad \forall \, \eta \in H^{1,1}_{0,per}(Q_T) \quad \forall \, z \in H^{1,\frac{1}{2}}_{0,per}(Q_T).
\end{align}
Also, we use the identities
\begin{align*}
 \int_\Omega \text{div} \, \boldsymbol{\rho} \, z \, d\boldsymbol{x} 
 = - \int_\Omega \boldsymbol{\rho} \cdot \nabla z \, d\boldsymbol{x} 
 \qquad \text{ and } \qquad
 \int_\Omega \text{div} \, \boldsymbol{\tau} \, q \, d\boldsymbol{x}
 = - \int_\Omega \boldsymbol{\tau} \cdot \nabla q \, d\boldsymbol{x}, 
\end{align*}
which hold for any $z,q \in H^1_0(\Omega)$ and any
\begin{align*}
 \boldsymbol{\tau}, \boldsymbol{\rho} \in H(\text{div},Q_T) :=
 \{\boldsymbol{\rho} \in [L^2(Q_T)]^d 
 : \text{div}_{\boldsymbol{x}} \, \boldsymbol{\rho}(\cdot,t) \in L^2(\Omega) 
 \text{ for a.e. } t \in (0,T)\}.
\end{align*}
For ease of notation, the index $\boldsymbol{x}$ in $\text{div}_{\boldsymbol{x}}$ will be henceforth omitted,
i.e., $\text{div} = \text{div}_{\boldsymbol{x}}$ denotes the generalized spatial divergence.
Using the Cauchy-Schwarz inequality yields the estimate 
\begin{align*}
 \mathcal{F}_{(\eta,\zeta)}(z,q&)
 = \int_{Q_T} \Big(y_d\,z
      - \eta\,z + \nu(\boldsymbol{x}) \nabla \zeta \cdot \nabla z
      - \sigma(\boldsymbol{x}) \partial_t \zeta \, z
      + \text{div} \, \boldsymbol{\rho} \, z + \boldsymbol{\rho} \cdot \nabla z \\
      &\qquad\, - \nu(\boldsymbol{x}) \nabla \eta \cdot \nabla q
      - \sigma(\boldsymbol{x}) \partial_t \eta \, q
      - \lambda^{-1} \zeta \,q
      + \text{div} \, \boldsymbol{\tau} \, q + \boldsymbol{\tau} \cdot \nabla q \Big) \, d\boldsymbol{x}\,dt \\
 \leq &\,\|\mathcal{R}_1(\eta,\zeta,\boldsymbol{\tau})\| 
 \|q\| + \|\mathcal{R}_2(\eta,\boldsymbol{\tau})\| \|\nabla q\| 
+ \|\mathcal{R}_3(\eta,\zeta,\boldsymbol{\rho})\| \|z\| 
 + \|\mathcal{R}_4(\zeta,\boldsymbol{\rho})\| \|\nabla z\|, 
 \end{align*}
where
\begin{align}
\label{definition:R1R2R3R4}
 \left.
\begin{aligned}
 \mathcal{R}_1(\eta,\zeta,\boldsymbol{\tau})
 &= \sigma \partial_t \eta + \lambda^{-1} \zeta - \text{div} \, \boldsymbol{\tau},
 \qquad \,\,\,
 \mathcal{R}_2(\eta,\boldsymbol{\tau}) = \boldsymbol{\tau}-\nu \nabla \eta, \\
  \mathcal{R}_3(\eta,\zeta,\boldsymbol{\rho})
 &= \sigma \partial_t \zeta + \eta - \text{div} \, \boldsymbol{\rho} - y_d,
 \qquad
 \mathcal{R}_4(\zeta,\boldsymbol{\rho}) = \boldsymbol{\rho} + \nu \nabla \zeta, 
\end{aligned}
\quad
\right \rbrace
\end{align}
Applying (\ref{inequality:Friedrichs:FourierSpace}), 
we find that
\begin{align*}
 \mathcal{F}_{(\eta,\zeta)}(z,q)
 &\leq  \Big(C_F \, \|\mathcal{R}_1(\eta,\zeta,\boldsymbol{\tau})\| 
 + \|\mathcal{R}_2(\eta,\boldsymbol{\tau})\| 
 \Big) \|\nabla q\| \\ 
 &+ \Big(C_F \, \|\mathcal{R}_3(\eta,\zeta,\boldsymbol{\rho})\| 
 + \|\mathcal{R}_4(\zeta,\boldsymbol{\rho})\| 
 \Big) \|\nabla z\|.  
\end{align*}
Hence,
\begin{align*}
  \sup_{0 \not= (z,q) \in (H^{1,\frac{1}{2}}_{0,per}(Q_T))^2}
    \frac{\mathcal{F}_{(\eta,\zeta)}(z,q)}{|(z,q)|_{1,\frac{1}{2}}} 
 \leq &\,  C_F \, \|\mathcal{R}_1(\eta,\zeta,\boldsymbol{\tau})\| 
 + \|\mathcal{R}_2(\eta,\boldsymbol{\tau})\|  \\ 
 &+ C_F \, \|\mathcal{R}_3(\eta,\zeta,\boldsymbol{\rho})\| 
 + \|\mathcal{R}_4(\zeta,\boldsymbol{\rho})\| 
\end{align*}
and by the 
inf-sup condition in (\ref{inequality:KKTSysSTVFinfsupsupsup:Seminorm}), we obtain
\begin{align*}
  |e|_{1,\frac{1}{2}} 
  \leq \frac{1}{\tilde \mu_1}
  \sup_{0 \not= (z,q) \in (H^{1,\frac{1}{2}}_{0,per}(Q_T))^2}
    \frac{\mathcal{B}(e, (z,q))}{|(z,q)|_{1,\frac{1}{2}}}
  = \frac{1}{\tilde \mu_1}
  \sup_{0 \not= (z,q) \in (H^{1,\frac{1}{2}}_{0,per}(Q_T))^2}
    \frac{\mathcal{F}_{(\eta,\zeta)}(z,q)}{|(z,q)|_{1,\frac{1}{2}}}.
\end{align*}
Thus, we arrive at the following result:
\begin{theorem}
\label{theorem:aposteriorEstimateH11/2Seminorm:OCP}
 Let $\eta,\zeta \in H^{1,1}_{0,per}(Q_T)$,
 $\boldsymbol{\tau}, \boldsymbol{\rho} \in H(\emph{div},Q_T)$,
  and the bilinear form $\mathcal{B}(\cdot,\cdot)$
 defined by (\ref{definition:KKTSysSTVF})
 meets the inf-sup condition (\ref{inequality:KKTSysSTVFinfsupsupsup:Seminorm}). Then,
 \begin{align}
 \left.
  \label{inequality:aposteriorEstimateH11/2Seminorm:OCP}
  \begin{aligned}
   |e|_{1,\frac{1}{2}} \leq \frac{1}{\tilde \mu_1}
   & \big(C_F \, \|\mathcal{R}_1(\eta,\zeta,\boldsymbol{\tau})\| 
     + \|\mathcal{R}_2(\eta,\boldsymbol{\tau})\|  \\
   &+ C_F \, \|\mathcal{R}_3(\eta,\zeta,\boldsymbol{\rho})\| 
     + \|\mathcal{R}_4(\zeta,\boldsymbol{\rho})\| 
     \big) 
   =: \mathcal{M}^\oplus_{|\cdot|}(\eta,\zeta,\boldsymbol{\tau},\boldsymbol{\rho}),
  \end{aligned}
\quad
\right \rbrace
 \end{align}
 where $e = (y,p) - (\eta, \zeta) \in (H^{1,\frac{1}{2}}_{0,per}(Q_T))^2$ and
 $\tilde \mu_1 = \frac{\min\{\underline{\nu},\underline{\sigma}\} \min\{\lambda,\frac{1}{\lambda}\}}{\sqrt{2}}$.
\end{theorem}
\begin{remark}
 For computational reasons, 
 it is useful to reformulate the majorants 
 in such a way that they are given by quadratic functionals,
 see, e.g., \cite{PhD:GaevskayaHoppeRepin:2006}. This is done by introducing 
 parameters 
 $\alpha, \beta, \gamma > 0$ and applying Young's inequality. 
 For the error majorant
 $\mathcal{M}^\oplus_{|\cdot|}(\eta,\zeta,\boldsymbol{\tau},\boldsymbol{\rho})$,
 we have 
 \begin{align*}
 \mathcal{M}^\oplus_{|\cdot|}(\eta,\zeta,\boldsymbol{\tau},\boldsymbol{\rho})^2
 \leq &\,\, \mathcal{M}^\oplus_{|\cdot|}(\alpha,\beta,\gamma;\eta,\zeta,\boldsymbol{\tau},\boldsymbol{\rho})^2 \\
  = &\, \frac{1}{\tilde \mu_1^2} \Big(C_F^2(1+\alpha)(1+\beta) \, \|\mathcal{R}_1(\eta,\zeta,\boldsymbol{\tau})\|^2 
     + \frac{(1+\alpha)(1+\beta)}{\beta} \|\mathcal{R}_2(\eta,\boldsymbol{\tau})\|^2   \\
   &+ C_F^2 \frac{(1+\alpha)(1+\gamma)}{\alpha} \, \|\mathcal{R}_3(\eta,\zeta,\boldsymbol{\rho})\|^2 
     + \frac{(1+\alpha)(1+\gamma)}{\alpha \gamma} \|\mathcal{R}_4(\zeta,\boldsymbol{\rho})\|^2 
     \Big). 
 \end{align*}
\end{remark}
Finally, we note that the inf-sup condition
 (\ref{inequality:KKTSysSTVFinfsupsupsup}) implies an estimate similar to
 (\ref{inequality:aposteriorEstimateH11/2Seminorm:OCP})
 for the error in terms of $\|\cdot\|_{1,\frac{1}{2}}$.
\begin{theorem}
\label{theorem:aposteriorEstimateH11/2Norm:OCP}
 Let $\eta,\zeta \in H^{1,1}_{0,per}(Q_T)$,
 $\boldsymbol{\tau},\boldsymbol{\rho} \in H(\emph{div},Q_T)$
  and the bilinear form $\mathcal{B}(\cdot,\cdot)$
 defined by (\ref{definition:KKTSysSTVF}) satisfies
 (\ref{inequality:KKTSysSTVFinfsupsupsup}). Then,
 \begin{align}
 \left.
 \label{inequality:aposteriorEstimateH11/2Nnorm:OCP}
 \begin{aligned}
  \|e\|_{1,\frac{1}{2}} \leq \frac{1}{\mu_{1}}
  &\Big(\|\mathcal{R}_1(\eta,\zeta,\boldsymbol{\tau})\|^2 
  + \|\mathcal{R}_2(\eta,\boldsymbol{\tau})\|^2  \\
  &+ \|\mathcal{R}_3(\eta,\zeta,\boldsymbol{\rho})\|^2
  + \|\mathcal{R}_4(\zeta,\boldsymbol{\rho})\|^2
  \Big)^{1/2}
  =: \mathcal{M}^\oplus_{\|\cdot\|}(\eta,\zeta,\boldsymbol{\tau},\boldsymbol{\rho}),
 \end{aligned}
\quad
\right \rbrace
 \end{align}
 where $e$ is defined in Theorem \ref{theorem:aposteriorEstimateH11/2Seminorm:OCP} and
 $\mu_1 = \frac{\min\{\frac{1}{\sqrt{\lambda}},\underline{\nu},\underline{\sigma}\}
 \min\{\sqrt{\lambda},\frac{1}{\sqrt{\lambda}}\}}{\sqrt{1+2\max\{\lambda,\frac{1}{\lambda}\}}}$.
\end{theorem}

\subsubsection*{The multiharmonic approximations}

Since the desired state $y_d$ belongs to $L^2$ 
it can be represented as a Fourier series. 
Henceforth, we assume that the
approximations $\eta$ and $\zeta$ of 
the exact state $y$ and the adjoint state $p$,
respectively,
are also represented in terms of truncated Fourier series
as well as the vector-valued functions
$\boldsymbol{\tau}$ and $\boldsymbol{\rho}$, i.e.,
\begin{align}
\label{definition:MhApproxEtaTau}
\left.
\begin{aligned}
 \eta(\boldsymbol{x},t) &= \eta_0^c(\boldsymbol{x}) + \sum_{k=1}^N \left(\eta_k^c(\boldsymbol{x}) \cos(k \omega t)
 + \eta_k^s(\boldsymbol{x}) \sin(k \omega t)\right), \\
 \boldsymbol{\tau}(\boldsymbol{x},t) &= \boldsymbol{\tau}_0^c(\boldsymbol{x}) + \sum_{k=1}^N \left(\boldsymbol{\tau}_k^c(\boldsymbol{x}) \cos(k \omega t)
 + \boldsymbol{\tau}_k^s(\boldsymbol{x}) \sin(k \omega t)\right),
\end{aligned}
\quad
\right \rbrace
\end{align}
where all the Fourier coefficients belong to the space $L^2(\Omega)$.
In this case, 
\begin{align*}
 \partial_t \eta(\boldsymbol{x},t) &= \sum_{k=1}^N \left(k \omega \, \eta_k^s(\boldsymbol{x}) \cos(k \omega t)
    - k \omega \, \eta_k^c(\boldsymbol{x}) \sin(k \omega t)\right), \\
 \nabla \eta(\boldsymbol{x},t) &= \nabla \eta_0^c(\boldsymbol{x})
 + \sum_{k=1}^N \left(\nabla \eta_k^c(\boldsymbol{x}) \, \cos(k \omega t)
 + \nabla \eta_k^s(\boldsymbol{x}) \, \sin(k \omega t)\right), \\
 \text{div} \, \boldsymbol{\tau}(\boldsymbol{x},t) &= \text{div} \, \boldsymbol{\tau}_0^c(\boldsymbol{x})
 + \sum_{k=1}^N \left(\text{div} \, \boldsymbol{\tau}_k^c(\boldsymbol{x}) \, \cos(k \omega t)
 + \text{div} \, \boldsymbol{\tau}_k^s(\boldsymbol{x}) \, \sin(k \omega t)\right),
\end{align*}
and the $L^2(Q_T)$-norms of the functions
$\mathcal{R}_1$, $\mathcal{R}_2$, $\mathcal{R}_3$  and $\mathcal{R}_4$ 
defined in (\ref{definition:R1R2R3R4}) can be 
represented in the form, which exposes each mode separately.
More precisely, we have
\begin{align*}
 \|\mathcal{R}_1(\eta,\zeta,\boldsymbol{\tau})\|^2
 = &\,T \|\lambda^{-1} \zeta_0^c - \text{div} \, \boldsymbol{\tau}_0^c\|_{\Omega}^2
 + \frac{T}{2} \sum_{k=1}^N
   \|-k \omega \, \sigma \boldsymbol{\eta}_k^\perp + \lambda^{-1} \boldsymbol{\zeta}_k 
      - \text{\textbf{div}} \, \boldsymbol{\tau}_k \|_{\Omega}^2, \\
 \|\mathcal{R}_2(\eta,\boldsymbol{\tau})\|^2 
 = &\,T \|\boldsymbol{\tau}_0^c - \nu \nabla \eta_0^c\|_{\Omega}^2 + \frac{T}{2} \sum_{k=1}^N
      \|\boldsymbol{\tau}_k - \nu \nabla \boldsymbol{\eta}_k\|_{\Omega}^2, \\
 \|\mathcal{R}_3(\eta,\zeta,\boldsymbol{\rho})\|^2
 = &\,T \|\eta_0^c - \text{div} \, \boldsymbol{\rho}_0^c - {y_d^c}_0\|_{\Omega}^2
 + \frac{T}{2} \sum_{k=1}^N
      \left(\|k \omega \, \sigma \zeta_k^s + \eta_k^c
       - \text{div} \, \boldsymbol{\rho}_k^c - {y_d^c}_k\|_{\Omega}^2 \right. \\
       &\left.+ \|-k \omega \, \sigma \zeta_k^c + \eta_k^s
       - \text{div} \, \boldsymbol{\rho}_k^s - {y_d^s}_k\|_{\Omega}^2\right)
 + \frac{T}{2} \sum_{k=N+1}^\infty
      \left(\|{y_d^c}_k\|_{\Omega}^2 + \|{y_d^s}_k\|_{\Omega}^2\right) \\
 = &\,T \|\eta_0^c - \text{div} \, \boldsymbol{\rho}_0^c - {y_d^c}_0\|_{\Omega}^2
 + \frac{T}{2} \sum_{k=1}^N
      \|-k \omega \, \sigma \boldsymbol{\zeta}_k^\perp + \boldsymbol{\eta}_k
      - \text{\textbf{div}} \, \boldsymbol{\rho}_k - \boldsymbol{y_d}_k\|_{\Omega}^2 \\
   &+ \frac{T}{2} \sum_{k=N+1}^\infty \|\boldsymbol{y_d}_k\|_{\Omega}^2, \\
 \|\mathcal{R}_4(\zeta,\boldsymbol{\rho})\|^2 
 = &\,T \|\boldsymbol{\rho}_0^c + \nu \nabla \zeta_0^c\|_{\Omega}^2 + \frac{T}{2} \sum_{k=1}^N
      \|\boldsymbol{\rho}_k + \nu \nabla \boldsymbol{\zeta}_k\|_{\Omega}^2,
\end{align*}
where $\text{\textbf{div}} \, \boldsymbol{\tau}_k
= (\text{div} \, \boldsymbol{\tau}_k^c,\text{div} \, \boldsymbol{\tau}_k^s)^T$ and
$\text{\textbf{div}} \, \boldsymbol{\rho}_k
= (\text{div} \, \boldsymbol{\rho}_k^c,\text{div} \, \boldsymbol{\rho}_k^s)^T$.
\begin{remark}
\label{remark:remainderterm}
Since $y_d$ is known, we can always compute the remainder term of truncation
\begin{align}
\label{def:remTerm}
\mathcal{E}_N :=
 \frac{T}{2} \sum_{k=N+1}^\infty \|\boldsymbol{y_d}_k\|_{\Omega}^2
 = \frac{T}{2} \sum_{k=N+1}^\infty \left(\|{y_d^c}_k\|_{\Omega}^2 + \|{y_d^s}_k\|_{\Omega}^2\right)
\end{align}
with any desired accuracy.
\end{remark}

It is important to outline that all the $L^2$-norms 
of $\mathcal{R}_1$, $\mathcal{R}_2$, $\mathcal{R}_3$ and
$\mathcal{R}_4$ corresponding to every single mode $k=0,\dots,N$ are decoupled.
It is useful to introduce quantities related to each mode. For $k \geq 1$, we denote them by
\begin{align*}
 {\mathcal{R}_1}_k(\boldsymbol{\eta}_k,\boldsymbol{\zeta}_k,\boldsymbol{\tau}_k)
 &:= -k \omega \, \sigma \boldsymbol{\eta}_k^\perp + \lambda^{-1} \boldsymbol{\zeta}_k 
      - \text{\textbf{div}} \, \boldsymbol{\tau}_k
 = ({\mathcal{R}_1}^c_k(\eta_k^s,\zeta_k^c,\boldsymbol{\tau}_k^c),
 {\mathcal{R}_1}^s_k(\eta_k^c,\zeta_k^s,\boldsymbol{\tau}_k^s))^T \\
 &= (k \omega \, \sigma \eta_k^s + \lambda^{-1} \zeta_k^c - \text{div} \, \boldsymbol{\tau}_k^c,
     -k \omega \, \sigma \eta_k^c + \lambda^{-1} \zeta_k^s - \text{div} \, \boldsymbol{\tau}_k^s)^T,
\end{align*}
\begin{align*}
 {\mathcal{R}_2}_k(\boldsymbol{\eta}_k,\boldsymbol{\tau}_k)
 &:= \boldsymbol{\tau}_k - \nu \nabla \boldsymbol{\eta}_k
 = ({\mathcal{R}_2}^c_k(\eta_k^c,\boldsymbol{\tau}_k^c),
 {\mathcal{R}_2}^s_k(\eta_k^s,\boldsymbol{\tau}_k^s))^T \\
 &= (\boldsymbol{\tau}_k^c - \nu \nabla \eta_k^c,
 \boldsymbol{\tau}_k^s - \nu \nabla \eta_k^s)^T, \\
 {\mathcal{R}_3}_k(\boldsymbol{\eta}_k,\boldsymbol{\zeta}_k,\boldsymbol{\rho}_k)
 &:= -k \omega \, \sigma \boldsymbol{\zeta}_k^\perp + \boldsymbol{\eta}_k
      - \text{\textbf{div}} \, \boldsymbol{\rho}_k - \boldsymbol{y_d}_k \\
 &= ({\mathcal{R}_3}^c_k(\eta_k^c,\zeta_k^s,\boldsymbol{\rho}_k^c),
 {\mathcal{R}_3}^s_k(\eta_k^s,\zeta_k^c,\boldsymbol{\rho}_k^s))^T \\
 &= (k \omega \, \sigma \zeta_k^s + \eta_k^c - \text{div} \, \boldsymbol{\rho}_k^c - {y_d^c}_k, 
 -k \omega \, \sigma \zeta_k^c + \eta_k^s - \text{div} \, \boldsymbol{\rho}_k^s - {y_d^s}_k)^T, \\
\text{and} \hspace{1.9 cm}& \\
 {\mathcal{R}_4}_k(\boldsymbol{\zeta}_k,\boldsymbol{\rho}_k)
 &:= \boldsymbol{\rho}_k + \nu \nabla \boldsymbol{\zeta}_k
 = ({\mathcal{R}_4}^c_k(\zeta_k^c,\boldsymbol{\rho}_k^c),
 {\mathcal{R}_4}^s_k(\zeta_k^s,\boldsymbol{\rho}_k^s))^T \\
 &= (\boldsymbol{\rho}_k^c + \nu \nabla \zeta_k^c,
    \boldsymbol{\rho}_k^s + \nu \nabla \zeta_k^s)^T.
\end{align*}
For $k=0$, we have 
\begin{align}
\label{definition:R1R2R3R4:Fourierk0}
\left.
 \begin{aligned}
 {\mathcal{R}_1}^c_0(\zeta_0^c,\boldsymbol{\tau}_0^c)
 &:= \lambda^{-1} \zeta_0^c - \text{div} \, \boldsymbol{\tau}_0^c, \qquad \qquad \quad
 {\mathcal{R}_2}^c_0(\eta_0^c,\boldsymbol{\tau}_0^c)
 := \boldsymbol{\tau}_0^c - \nu \nabla \eta_0^c, \\  
 {\mathcal{R}_3}^c_0(\eta_0^c,\boldsymbol{\rho}_0^c)
 &:= \eta_0^c - \text{div} \, \boldsymbol{\rho}_0^c - {y_d^c}_0, \qquad \qquad
 {\mathcal{R}_4}^c_0(\zeta_0^c,\boldsymbol{\rho}_0^c)
 := \boldsymbol{\rho}_0^c + \nu \nabla \zeta_0^c.
 \end{aligned}
 \right \rbrace
\end{align}
\begin{corollary}
The error majorants $\mathcal{M}^\oplus_{|\cdot|}(\eta,\zeta,\boldsymbol{\tau},\boldsymbol{\rho})$ and
$\mathcal{M}^\oplus_{\|\cdot\|}(\eta,\zeta,\boldsymbol{\tau},\boldsymbol{\rho})$ 
defined in (\ref{inequality:aposteriorEstimateH11/2Seminorm:OCP}) and
(\ref{inequality:aposteriorEstimateH11/2Nnorm:OCP}),
respectively,
can be represented 
in somewhat new forms that contain quantities associated with the modes, namely, 
 \begin{align*}
  \begin{aligned}
  \mathcal{M}^\oplus_{|\cdot|}(\eta,\zeta,\boldsymbol{\tau},\boldsymbol{\rho})
  = \frac{1}{\tilde \mu_1} 
   &\Big(C_F \, \big(T \|{\mathcal{R}_1}^c_0(\zeta_0^c,\boldsymbol{\tau}_0^c)\|_{\Omega}^2
  + \frac{T}{2} \sum_{k=1}^N
   \|{\mathcal{R}_1}_k(\boldsymbol{\eta}_k,\boldsymbol{\zeta}_k,\boldsymbol{\tau}_k)\|_{\Omega}^2 \big)^{1/2} \\
     &+ \big( T \|{\mathcal{R}_2}^c_0(\eta_0^c,\boldsymbol{\tau}_0^c)\|_{\Omega}^2
     + \frac{T}{2} \sum_{k=1}^N \|{\mathcal{R}_2}_k(\boldsymbol{\eta}_k,\boldsymbol{\tau}_k)\|_{\Omega}^2
     \big)^{1/2} \\
   &+ \,C_F \, \big(T \|{\mathcal{R}_3}^c_0(\eta_0^c,\boldsymbol{\rho}_0^c)\|_{\Omega}^2
   + \frac{T}{2} \sum_{k=1}^N
      \|{\mathcal{R}_3}_k(\boldsymbol{\eta}_k,\boldsymbol{\zeta}_k,\boldsymbol{\rho}_k)\|_{\Omega}^2
   + \mathcal{E}_N
     \big)^{1/2} \\
   &+\, \big(T \|{\mathcal{R}_4}^c_0(\zeta_0^c,\boldsymbol{\rho}_0^c)\|_{\Omega}^2 + \frac{T}{2} \sum_{k=1}^N
      \|{\mathcal{R}_4}_k(\boldsymbol{\zeta}_k,\boldsymbol{\rho}_k)\|_{\Omega}^2\big)^{1/2} \Big) 
  \end{aligned}
 \end{align*}
 and
 \begin{align*}
 \begin{aligned}
 \mathcal{M}^\oplus_{\|\cdot\|}(\eta,\zeta,\boldsymbol{\rho}&,\boldsymbol{\tau})
  = \frac{1}{\mu_{1}}
  \Big(T \big(\|{\mathcal{R}_1}^c_0(\zeta_0^c,\boldsymbol{\tau}_0^c)\|_{\Omega}^2 
  + \|{\mathcal{R}_2}^c_0(\eta_0^c,\boldsymbol{\tau}_0^c)\|_{\Omega}^2
  + \|{\mathcal{R}_3}^c_0(\eta_0^c,\boldsymbol{\rho}_0^c)\|_{\Omega}^2
   \\
  &+ \|{\mathcal{R}_4}^c_0(\zeta_0^c,\boldsymbol{\rho}_0^c)\|_{\Omega}^2 \big)
  + \frac{T}{2} \sum_{k=1}^N
  \big(\|{\mathcal{R}_1}_k(\boldsymbol{\eta}_k,\boldsymbol{\zeta}_k,\boldsymbol{\tau}_k)\|_{\Omega}^2 
  + \|{\mathcal{R}_2}_k(\boldsymbol{\eta}_k,\boldsymbol{\tau}_k)\|_{\Omega}^2 \\
  &+ \|{\mathcal{R}_3}_k(\boldsymbol{\eta}_k,\boldsymbol{\zeta}_k,\boldsymbol{\rho}_k)\|_{\Omega}^2 
  + \|{\mathcal{R}_4}_k(\boldsymbol{\zeta}_k,\boldsymbol{\rho}_k)\|_{\Omega}^2\big)
  + \mathcal{E}_N \Big)^{1/2}. 
 \end{aligned}
 \end{align*}
\end{corollary}
\begin{remark}
Let $y_d$ has a multiharmonic representation, i.e.,
 \begin{align}
 \label{def:ydMultiharmonic}
  y_d(\boldsymbol{x},t) = {y_d}_0^c(\boldsymbol{x}) + \sum_{k=1}^{N_d} \left(f_k^c(\boldsymbol{x}) \cos(k \omega t)
  + {y_d}_k^s(\boldsymbol{x}) \sin(k \omega t)\right),
 \end{align}
where $N_d \in \mathbb{N}$.
If $N \geq N_d$, then 
$(\eta,\zeta)$ is the exact solution of problem (\ref{problem:KKTSysSTVFAPost}) and
$(\boldsymbol{\tau},\boldsymbol{\rho})$ is the exact flux if and only if the error majorants 
vanish, i.e.,
 \begin{align*}
  {\mathcal{R}_j}_0^c = 0 \qquad \text{and} \qquad
  {\mathcal{R}_j}_k = 0 \qquad \forall \, k=1,\dots,N_d, \qquad \forall \, j \in \{1,2,3,4\}.
 \end{align*}
Indeed, let the error majorants vanish. Then,
$\eta_0^c - \text{\emph{div}} \, \boldsymbol{\rho}_0^c = {y_d^c}_0$,  
$\boldsymbol{\rho}_0^c = - \nu \nabla \zeta_0^c$,  
$\lambda^{-1} \zeta_0^c - \text{\emph{div}} \, \boldsymbol{\tau}_0^c = 0$,
$\boldsymbol{\tau}_0^c = \nu \nabla \eta_0^c$
and we see that
 \begin{align*}
  &k \omega \, \sigma \eta_k^s + \lambda^{-1} \zeta_k^c - \text{\emph{div}} \, \boldsymbol{\tau}_k^c = 0, \qquad
     -k \omega \, \sigma \eta_k^c + \lambda^{-1} \zeta_k^s - \text{\emph{div}} \, \boldsymbol{\tau}_k^s = 0, \\
  &k \omega \, \sigma \zeta_k^s + \eta_k^c - \text{\emph{div}} \, \boldsymbol{\rho}_k^c = {y_d^c}_k,
  \qquad \quad \, \,
   -k \omega \, \sigma \zeta_k^c + \eta_k^s - \text{\emph{div}} \, \boldsymbol{\rho}_k^s = {y_d^s}_k, \\
    &\boldsymbol{\tau}_k^c = \nu \nabla \eta_k^c, \qquad
    \boldsymbol{\tau}_k^s = \nu \nabla \eta_k^s, \qquad
    \boldsymbol{\rho}_k^c = - \nu \nabla \zeta_k^c, \qquad
    \boldsymbol{\rho}_k^s = - \nu \nabla \zeta_k^s, 
 \end{align*}
for all $k=1,\dots,N_d$, 
so that collecting the $N_d+1$ harmonics leads to
multiharmonic representations for $\eta$, $\zeta$, $\boldsymbol{\tau}$
and $\boldsymbol{\rho}$ of the form (\ref{def:ydMultiharmonic}) 
satisfying the equations 
\begin{align*}
  \sigma \partial_t \eta - \emph{div} \, \boldsymbol{\tau} + \lambda^{-1} \zeta = 0, \qquad
   \boldsymbol{\tau} = \nu \nabla \eta, \qquad
  \sigma \partial_t \zeta - \emph{div} \, \boldsymbol{\rho} + \eta = y_d, \qquad
   \boldsymbol{\rho} = - \nu \nabla \zeta.
\end{align*}
Since $\eta$ and $\zeta$ also meet the boundary conditions,
we conclude that $\eta = y$ and $\zeta = p$.
\end{remark}

\section[A posteriori estimates for cost functionals]{Functional
A Posteriori Estimates for Cost Functionals of Parabolic Time-Periodic Optimal Control Problems}
\label{Sec5:APosterioriErrorEstimation:OCP:CostFuncs}

This section is aimed at deriving guaranteed and computable
upper bounds for the cost functional
 and establishing their sharpness.
 This is important because, in optimal control, we cannot in general compute the (exact) 
 cost functional since
 the (exact) state function is unknown. By using {\it a posteriori} estimates for the state equation
 we overcome this difficulty.
 Similar results for elliptic optimal control problems can be found, e.g., in
 \cite{PhD:GaevskayaHoppeRepin:2006, PhD:Repin:2008}.
 Let $y=y(v)$ be the corresponding state to a control $v$.
 The cost functional $\mathcal{J}(y(v),v)$ 
 defined in (\ref{equation:minfunc:OCP}) has the form 
 \begin{align*}
  \mathcal{J}(y(v),v) 
	=  T \mathcal{J}_0(y_0^c(v_0^c),v_0^c)
	+ \frac{T}{2} \sum_{k=1}^\infty \mathcal{J}_k(\boldsymbol{y}_k(\boldsymbol{v}_k),\boldsymbol{v}_k),
 \end{align*}
 where $\mathcal{J}_0(y_0^c(v_0^c),v_0^c)
  = \frac{1}{2} \|y_0^c - {y_d^c}_0\|_{\Omega}^2
     + \frac{\lambda}{2} \|v_0^c\|_{\Omega}^2$
 and
 \begin{align*}
 \mathcal{J}_k(\boldsymbol{y}_k(\boldsymbol{v}_k),\boldsymbol{v}_k)
 = \frac{1}{2} \|\boldsymbol{y}_k-{\boldsymbol{y}_d}_k\|_{\Omega}^2
     + \frac{\lambda}{2} \|\boldsymbol{v}_k\|_{\Omega}^2.
 \end{align*}
 We wish to deduce majorants for the cost functional
  $\mathcal{J}(y(u),u)$ of the exact control $u$ and corresponding state $y(u)$
  by using some of the results presented in \cite{LRW:LangerRepinWolfmayr:2015},
  which are obtained for the time-periodic boundary value problem (\ref{equation:forwardpde:OCP}).
In \cite{LRW:LangerRepinWolfmayr:2015},
the following functional {\it a posteriori} error estimate for problem
(\ref{equation:forwardpde:OCP}) has been proved:
 \begin{align}
  \label{inequality:aposteriorEstimateH11/2SeminormBVP}
  |y(v) - \eta|_{1,\frac{1}{2}} \leq \frac{1}{\underline{\mu_1}}
  \left(C_F \, \|\mathcal{R}_1(\eta,\boldsymbol{\tau},v)\| 
    + \|\mathcal{R}_2(\eta,\boldsymbol{\tau})\| \right),
 \end{align}
 where
 $\underline{\mu_{1}} = \frac{1}{\sqrt{2}} \min\{\underline{\nu},\underline{\sigma}\}$.
 It holds for arbitrary functions $\eta \in H^{1,1}_{0,per}(Q_T)$ and
 $\boldsymbol{\tau} \in H(\text{div},Q_T)$,
 where
 \begin{align*}
 \mathcal{R}_1(\eta,\boldsymbol{\tau},v)
   := \sigma \partial_t \eta - \text{div} \, \boldsymbol{\tau} - v,  \qquad
   \qquad
 \mathcal{R}_2(\eta,\boldsymbol{\tau})
   := \boldsymbol{\tau} - \nu \nabla \eta,
  \end{align*}
and $v$ is a given function in $L^2(Q_T)$.
Now, adding and subtracting $\eta$ in the cost functional $\mathcal{J}(y(v),v)$ as well as
  applying the triangle and Friedrichs inequalities yields the estimate 
  \begin{align*}
   \mathcal{J}(y(v),v) 
   &\leq \frac{1}{2} \left(\|\eta - y_d\| + C_F \|\nabla y(v) - \nabla \eta\| \right)^2
     + \frac{\lambda}{2} \|v\|^2.
  \end{align*}
  Since
  \begin{align*}
   \|\nabla y(v) - \nabla \eta\|^2
    \leq |y(v) - \eta|_{1,\frac{1}{2}}^2 = \|\nabla y(v) - \nabla \eta\|^2
    + \|\partial_t^{1/2} y(v) - \partial_t^{1/2} \eta\|^2,
  \end{align*}
  we conclude that 
  \begin{align*}
   \mathcal{J}(y(v),v)
   &\leq \frac{1}{2} \left(\|\eta - y_d\| + C_F |y(v) - \eta|_{1,\frac{1}{2}} \right)^2
     + \frac{\lambda}{2} \|v\|^2.
  \end{align*}  
  Together with (\ref{inequality:aposteriorEstimateH11/2SeminormBVP}) 
this leads to the estimate 
  \begin{align*}
   \mathcal{J}(y(v),v)
   \leq \frac{1}{2} \left(\|\eta - y_d\|
     + \frac{C_F}{\underline{\mu_1}} \|\mathcal{R}_2(\eta,\boldsymbol{\tau})\|
     + \frac{C_F^2}{\underline{\mu_1}} \|\mathcal{R}_1(\eta,\boldsymbol{\tau},v)\| \right)^2
     + \frac{\lambda}{2} \|v\|^2.
  \end{align*}
  By introducing parameters $\alpha, \beta > 0$ and applying Young's inequality,
  we can reformulate the estimate such that the
  right-hand side is given by a quadratic functional 
  (the latter functional is more convenient from the computational point of view).
We have
  \begin{align*}
   \mathcal{J}(y(v),v)
   \leq \mathcal{J}^\oplus(\alpha,\beta;\eta,\boldsymbol{\tau},v) \qquad \forall \, v \in L^2(Q_T),
  \end{align*}
  where
  \begin{align}
  \left.
  \label{definition:majorantCostFunc}
  \begin{aligned}
   \mathcal{J}^\oplus(\alpha,\beta;\eta,\boldsymbol{\tau},v) :=
   &\, \frac{1+\alpha}{2} \|\eta - y_d\|^2
     + \frac{(1+\alpha)(1+\beta) C_F^2 }{2 \alpha \underline{\mu_1}^2} \|\mathcal{R}_2(\eta,\boldsymbol{\tau})\|^2 \\
     &+ \frac{(1+\alpha)(1+\beta) C_F^4}{2 \alpha \beta \underline{\mu_1}^2}
      \|\mathcal{R}_1(\eta,\boldsymbol{\tau},v)\|^2 
     + \frac{\lambda}{2} \|v\|^2.
  \end{aligned}
\quad
\right \rbrace
  \end{align}
The majorant (\ref{definition:majorantCostFunc}) provides a guaranteed 
upper bound of the cost functional,
which can be computed for any approximate control and state functions.
Moreover, minimization of this functional with respect to
$\eta$, $\boldsymbol{\tau}$, $v$ and $\alpha, \beta > 0$
yields the exact value of the cost functional.
This important result is summarized in the following theorem:
\begin{theorem}
 \label{theorem:majorantCostFuncInf}
 The exact lower bound of the majorant $\mathcal{J}^\oplus$ defined in (\ref{definition:majorantCostFunc}) coincides 
 with the optimal value of the cost functional of problem (\ref{equation:minfunc:OCP})-(\ref{equation:forwardpde:OCP}),
 i.e.,
 \begin{align}
  \label{definition:majorantCostFuncInf}
   \inf_{\substack{\eta \in H^{1,1}_{0,per}(Q_T),\boldsymbol{\tau} \in H(\emph{div},Q_T) \\
   v \in L^2(Q_T), \alpha, \beta > 0}}
   \mathcal{J}^\oplus(\alpha,\beta;\eta,\boldsymbol{\tau},v) = \mathcal{J}(y(u),u).
 \end{align}
\begin{proof}
 The infimum of $\mathcal{J}^\oplus$ is attained for the optimal control $u$, the 
 corresponding state function $y(u)$ and the exact flux $(\nu \nabla y(u))$.
  In this case, $\mathcal{R}_1$ and $\mathcal{R}_2$ vanish, and for $\alpha$ tending
  to zero, values of $\mathcal{J}^\oplus$ tend to the exact value of the
  cost functional.
\end{proof}
\end{theorem}
\begin{corollary}
From Theorem~\ref{theorem:majorantCostFuncInf}, we obtain
  the following estimate:
  \begin{align}
  \label{estimate:majCostFunc}
  \begin{aligned}
   \mathcal{J}(y(u),u)
   \leq \,\, &\mathcal{J}^\oplus(\alpha,\beta;\eta,\boldsymbol{\tau},v) \\ 
   &\forall \, \eta \in H^{1,1}_{0,per}(Q_T), \, 
   \boldsymbol{\tau} \in H(\emph{div},Q_T), \, 
   v \in L^2(Q_T), \, \alpha, \beta > 0.
  \end{aligned}
  \end{align}
\end{corollary}
  Now, it is easy to derive {\it a posteriori} 
  estimates for the cost functional in the setting of multiharmonic approximations.
  Let $\eta$ be the MhFE approximation $y_{N h}$ to the state $y$. 
  Since the control $v$ can be chosen arbitrarily in (\ref{definition:majorantCostFunc}), 
  we choose a MhFE approximation
  $u_{N h}$ for it as well. More precisely, 
  we can compute the MhFE approximation
  $u_{N h}$ for the control from the MhFE approximation
  $p_{N h}$ of the adjoint state, since $u_{N h} = - \lambda^{-1} p_{N h}$, 
  by solving the optimality system,
  from which we obtain $y_{N h}$ as well.
We now apply (\ref{estimate:majCostFunc}) and select $\eta = y_{N h}$ and
$v = u_{N h}$.
  Next, we need to make a suitable reconstruction of 
  $\boldsymbol{\tau}$, which can be done by different techniques,
  see, e.g., \cite{PhD:Repin:2008, PhD:MaliNeittaanmaekiRepin:2014} and the references therein.
In the treatment of our approach it is natural to represent $\boldsymbol{\tau}$ in the form of
a multiharmonic function $\boldsymbol{\tau}_{N h}$. Then the majorant 
  \begin{align*}
  \mathcal{J}^\oplus&(\alpha,\beta;y_{N h},
  \boldsymbol{\tau}_{N h},u_{N h}) =
     \frac{1+\alpha}{2} \|y_{N h} - y_d\|^2
     + \frac{(1+\alpha)(1+\beta) C_F^2 }{2 \alpha \underline{\mu_1}^2} \|\mathcal{R}_2(y_{N h},
     \boldsymbol{\tau}_{N h})\|^2 \\
     &\qquad \qquad \qquad \qquad \qquad \, + \frac{(1+\alpha)(1+\beta) C_F^4}{2 \alpha 
     \beta \underline{\mu_1}^2}
     \|\mathcal{R}_1(y_{N h},
     \boldsymbol{\tau}_{N h},u_{N h})\|^2 
     + \frac{\lambda}{2} \|u_{N h}\|^2
  \end{align*}    
  has a multiharmonic form
  \begin{align}
  \label{def:costfunc:multiharmonic}
  \begin{aligned}
  &\,\mathcal{J}^\oplus(\alpha,\beta;y_{N h},
  \boldsymbol{\tau}_{N h},u_{N h})
  = \, \frac{1+\alpha}{2} \Big(T \|y_{0 h}^c
  - {y_d}_0^c\|_{\Omega}^2 \\
      &\qquad \qquad \qquad \qquad \quad 
      + \frac{T}{2} \sum_{k=1}^N
      \left(\|y_{k h}^c - {y_d}_k^c\|_{\Omega}^2
      + \|y_{k h}^s - {y_d}_k^s\|_{\Omega}^2 \right) + \mathcal{E}_N
      \Big) \\
     &+ \frac{(1+\alpha)(1+\beta) C_F^2 }{2 \alpha \underline{\mu_1}^2}
     \Big(T \|{\mathcal{R}_2}^c_0(y_{0 h}^c,
     \boldsymbol{\tau}_{0 h}^c)\|_{\Omega}^2 
     \\
     &\qquad \qquad \qquad \qquad \quad 
     + \frac{T}{2} \sum_{k=1}^N
      \left(\|{\mathcal{R}_2}^c_k(y_{k h}^c,
      \boldsymbol{\tau}_{k h}^c)\|_{\Omega}^2
      + \|{\mathcal{R}_2}^s_k(y_{k h}^s,
      \boldsymbol{\tau}_{k h}^s)\|_{\Omega}^2\right)\Big) \\
     &+ \frac{(1+\alpha)(1+\beta) C_F^4}{2 \alpha \beta \underline{\mu_1}^2}
     \Big(T \|{\mathcal{R}_1}^c_0(\boldsymbol{\tau}_{0 h}^c,u_{0 h}^c)\|_{\Omega}^2 \\
      &\qquad \qquad \qquad 
      + \frac{T}{2} \sum_{k=1}^N
      \left(\|{\mathcal{R}_1}^c_k(y_{k h}^s,
      \boldsymbol{\tau}_{k h}^c,u_{k h}^c)\|_{\Omega}^2
      + \|{\mathcal{R}_1}^s_k(y_{k h}^c,
      \boldsymbol{\tau}_{k h}^s,u_{k h}^s)\|_{\Omega}^2\right)\Big) \\
     &+ \frac{\lambda}{2} \Big(T \|u_{0 h}^c\|_{\Omega}^2
     + \frac{T}{2} \sum_{k=1}^N 
      \left(\|u_{k h}^c\|_{\Omega}^2
      + \|u_{k h}^s\|_{\Omega}^2\right)\Big),
 \end{aligned}
\end{align}
where
${\mathcal{R}_1}^c_0(\boldsymbol{\tau}_{0 h}^c,u_{0 h}^c) = \text{div} \, \boldsymbol{\tau}_{0 h}^c + u_{0 h}^c$,
${\mathcal{R}_2}^c_0(y_{0 h}^c,\boldsymbol{\tau}_{0 h}^c) = \boldsymbol{\tau}_{0 h}^c - \nu \nabla y_{0 h}^c$,
\begin{align*}
 {\mathcal{R}_1}_k(\boldsymbol{y}_{k h},
 \boldsymbol{\tau}_{k h},\boldsymbol{u}_{k h})
  &= k \omega \, \sigma \boldsymbol{y}_{k h}^\perp
  + \text{\textbf{div}} \, \boldsymbol{\tau}_{k h} + \boldsymbol{u}_{k h} \\
  &= (-k \omega \, \sigma y_{k h}^s 
  + \text{div} \, \boldsymbol{\tau}_{k h}^c + u_{k h}^c,
     k \omega \, \sigma y_{k h}^c 
     + \text{div} \, \boldsymbol{\tau}_{k h}^s + u_{k h}^s)^T \\
  &= ({\mathcal{R}_1}^c_k(y_{k h}^s,
      \boldsymbol{\tau}_{k h}^c,u_{k h}^c),
      {\mathcal{R}_1}^s_k(y_{k h}^c,
      \boldsymbol{\tau}_{k h}^s,u_{k h}^s))^T,
\end{align*}  
and
\begin{align*}
 {\mathcal{R}_2}_k(\boldsymbol{y}_{k h},\boldsymbol{\tau}_{k h})
  &= \boldsymbol{\tau}_{k h} - \nu \nabla \boldsymbol{y}_{k h} \\ 
  &= (\boldsymbol{\tau}_{k h}^c - \nu \nabla y_{k h}^c, 
     \boldsymbol{\tau}_{k h}^s - \nu \nabla y_{k h}^s)^T 
     \\
  &= ({\mathcal{R}_2}^c_k(y_{k h}^c,
      \boldsymbol{\tau}_{k h}^c),
      {\mathcal{R}_2}^s_k(y_{k h}^s,
      \boldsymbol{\tau}_{k h}^s))^T.
\end{align*}
Note that the 
remainder term (\ref{def:remTerm}) remains 
the same in (\ref{def:costfunc:multiharmonic}).
 Since all the terms corresponding to every single mode $k$ in the majorant $\mathcal{J}^\oplus$ are decoupled,
 we arrive at some majorants $\mathcal{J}_k^\oplus$, for which 
 we can, of course, introduce positive parameters $\alpha_k$ and
 $\beta_k$ for every single mode $k$ as well. Then the majorant (\ref{def:costfunc:multiharmonic}) 
 can be written as
   \begin{align}
  \label{def:costfunc:multiharmonic:new}
  \begin{aligned}
  \mathcal{J}^\oplus(\boldsymbol{\alpha}_{N+1},\boldsymbol{\beta}_{N};y_{N h},
  \boldsymbol{\tau}_{N h},u_{N h})
  = &\, T \,  \mathcal{J}_0^\oplus(\alpha_0,\beta_0;y_{0 h}^c,\boldsymbol{\tau}_{0 h}^c,u_{0 h}^c) \\
      &+ \frac{T}{2} \sum_{k=1}^N 
 \mathcal{J}_k^\oplus(\alpha_k,\beta_k; \boldsymbol{y}_{k h}, \boldsymbol{\tau}_{k h}, \boldsymbol{u}_{k h}) \\
&+ \frac{1+\alpha_{N+1}}{2} \, \mathcal{E}_N,
 \end{aligned}
\end{align}
where $\boldsymbol{\alpha}_{N+1} = (\alpha_0,\dots,\alpha_{N+1})^T$,
$\boldsymbol{\beta}_N = (\beta_0,\dots,\beta_N)^T$, and
   \begin{align}
    \begin{aligned}
    \label{def:Joplus0}
  \mathcal{J}_0^\oplus(\alpha_0,\beta_0&;y_{0 h}^c,
  \boldsymbol{\tau}_{0 h}^c,u_{0 h}^c)
  = \frac{1+\alpha_0}{2} \|y_{0 h}^c
  - {y_d}_0^c\|_{\Omega}^2 + \frac{\lambda}{2} \|u_{0 h}^c\|_{\Omega}^2 \\
  &+ \frac{(1+\alpha_0)(1+\beta_0) C_F^2 }{2 \alpha_0 \underline{\mu_1}^2}
      \|{\mathcal{R}_2}^c_0 \|_{\Omega}^2
  + \frac{(1+\alpha_0)(1+\beta_0) C_F^4}{2 \alpha_0 \beta_0 \underline{\mu_1}^2}
     \|{\mathcal{R}_1}^c_0 \|_{\Omega}^2,
\end{aligned}
\end{align}
   \begin{align}
    \begin{aligned}      
    \label{def:Joplusk}
  \mathcal{J}_k^\oplus(\alpha_k,\beta_k&; \boldsymbol{y}_{k h},
  \boldsymbol{\tau}_{k h}, \boldsymbol{u}_{k h})
  = \frac{1+\alpha_k}{2} \|\boldsymbol{y}_{k h} 
      - {\boldsymbol{y}_d}_k \|_{\Omega}^2
	+ \frac{\lambda}{2}   \|\boldsymbol{u}_{k h} \|_{\Omega}^2 \\
      &+ \frac{(1+\alpha_k)(1+\beta_k) C_F^2 }{2 \alpha_k \underline{\mu_1}^2}
      \|{\mathcal{R}_2}_k \|_{\Omega}^2
      + \frac{(1+\alpha_k)(1+\beta_k) C_F^4}{2 \alpha_k \beta_k \underline{\mu_1}^2}
      \|{\mathcal{R}_1}_k \|_{\Omega}^2. 
 \end{aligned}
\end{align}
Next, we have to reconstruct the fluxes $\boldsymbol{\tau}_{0 h}^c$ and
$\boldsymbol{\tau}_{k h}$ for all $k=1,\dots,N$, which we denote by
\begin{align*}
 \boldsymbol{\tau}_{k h}
 = R_h^{\text{\tiny{flux}}}(\nu \nabla \boldsymbol{y}_{k h}).
\end{align*}
This can be done by various techniques.
In \cite{LRW:LangerRepinWolfmayr:2015},
we have used Raviart-Thomas elements of the lowest order 
(see also 
\cite{PhD:RaviartThomas:1977, LRW:BrezziFortin:1991, LRW:RobertsThomas:1991}),
in order to regularize the fluxes by a post-processing operator, which maps the $L^2$-functions
into $H(\text{div},\Omega)$.
Collecting all the fluxes corresponding to the 
modes together yields
the reconstructed flux
\begin{align*}
 \boldsymbol{\tau}_{N h} 
 = R_h^{\text{\tiny{flux}}}(\nu \nabla y_{N h}).
\end{align*}
After performing a simple minimization of the majorant
$\mathcal{J}^\oplus$
with respect to 
$\boldsymbol{\alpha}_{N+1}$ and $\boldsymbol{\beta}_{N}$,
we finally arrive at the {\it a posteriori} estimate
  \begin{align}
   \begin{aligned}
   \mathcal{J}(y(u),u)
   \leq \mathcal{J}^\oplus(\boldsymbol{\bar{\alpha}}_{N+1},\boldsymbol{\bar{\beta}}_{N};
   y_{N h},\boldsymbol{\tau}_{N h},u_{N h}),
  \end{aligned}
  \end{align}
where $\boldsymbol{\bar{\alpha}}_{N+1}$ and $\boldsymbol{\bar{\beta}}_{N}$
denote the optimized positive parameters.

It is worth outlining that the 
majorant $\mathcal{J}^{\oplus}$ provides a guaranteed upper bound for the cost functional,
and, due to Theorem \ref{theorem:majorantCostFuncInf}, the infimum of the majorant
coincides with the optimal value of the cost functional.
\begin{remark}
 In this work,
 problems with constraints on the control or the state are not considered,
 but inequality constraints imposed on the Fourier coefficients of the control
 can easily be included into the MhFE approach, see
 \cite{PhD:KollmannKolmbauer:2011},
 and, hence, may be considered in the {\it a posteriori} error analysis of parabolic 
 time-periodic optimal control problems as well.
\end{remark}

\section{Numerical Results}
\label{Sec6:NumericalResults}

We compute and analyze the efficiency of the above derived
{\it a posteriori} estimates for different cases, namely,
\begin{enumerate}
 \item[{1.}] the desired state is periodic and analytic in time, but not time-harmonic,
 \item[{2.}] the desired state is analytic in time, but not time-periodic, and
 \item[{3.}] the desired state is a non-smooth function in space and time. 
\end{enumerate}
Note that convergence and other properties of numerical approximations generated
by the MhFEM have been studied in 
\cite{LRW:KollmannKolmbauerLangerWolfmayrZulehner:2013, LRW:LangerWolfmayr:2013}
for the same three cases.
The optimal control problem 
(\ref{equation:minfunc:OCP})-(\ref{equation:forwardpde:OCP}) is solved 
on the computational domain 
$\Omega = (0,1) \times (0,1)$
with the Friedrichs constant $C_F = 1/(\sqrt{2}\pi)$
using a uniform simplicial mesh 
and standard continuous, piecewise linear finite elements. 
The material coefficients are supposed to be 
$\sigma = \nu = 1$.
In the first two examples, we choose the cost parameter $\lambda = 0.1$, and
in the third one, we choose $\lambda = 0.01$. Both choices are common.

Getting sharp error bounds requires an efficient construction of
$\eta$, $\zeta$, $\boldsymbol{\tau}$ and  $\boldsymbol{\rho}$
in order to compute sharp guaranteed bounds from the majorants.
We choose MhFE approximations
(\ref{definition:MhApproxEtaTau}) for $\eta$ and $\boldsymbol{\tau}$ as well as for
$\zeta$ and $\boldsymbol{\rho}$.
In order to obtain
suitable fluxes $\boldsymbol{\tau}, \boldsymbol{\rho} \in H(\text{div},Q_T)$,
we reconstruct them 
by the standard lowest-order Raviart-Thomas ($RT^0$-) extension of normal fluxes.
We refer the reader to \cite{LRW:LangerRepinWolfmayr:2015}, where the authors
have discussed this issue thoroughly.

In order to solve the saddle point systems (\ref{equation:MultiFESysBlockk}) for $k=1,\dots,N$
and (\ref{equation:MultiFESysBlock0}) for $k=0$,
we use the robust algebraic multilevel  preconditioner of Kraus (see \cite{LRW:Kraus:2012})
for an inexact realization of the block-diagonal preconditioners
\begin{align}
\label{equation:preconditionerSigmaPWConst:OCP}
 \mathcal{P}_k = \left( \begin{array}{cccc}
     D & 0 & 0 & 0 \\
     0 & D & 0 & 0 \\
     0 & 0 & \lambda^{-1} D & 0 \\
     0 & 0 & 0 & \lambda^{-1} D \end{array} \right),
\end{align}
where 
$D = \sqrt{\lambda} K_{h,\nu} + k \omega \sqrt{\lambda} M_{h,\sigma} + M_h,$
and
\begin{align}
\label{equation:preconditionerSigmaPWConstkIs0:OCP}
 \mathcal{P}_0 = \left( \begin{array}{cc}
     M_h + \sqrt{\lambda} K_{h,\nu} & 0 \\
     0 & \lambda^{-1} (M_h + \sqrt{\lambda} K_{h,\nu}) \end{array} \right),
\end{align}
in the minimal residual method, respectively. 
The preconditioners (\ref{equation:preconditionerSigmaPWConst:OCP}) and
(\ref{equation:preconditionerSigmaPWConstkIs0:OCP})
were presented and discussed in
\cite{LRW:KollmannKolmbauerLangerWolfmayrZulehner:2013, LRW:LangerWolfmayr:2013, 
LRW:Wolfmayr:2014}.
The numerical experiments were computed on grids of  
mesh sizes 16$\times$16 to 256$\times$256.
The algorithms were implemented in C\texttt{++}, and
all computations were performed on an average class 
computer with Intel(R) Xeon(R) CPU W3680 @ 3.33 GHz.
Note that the 
presented CPU times $t^{\text{sec}}$ in seconds 
include the computational times for computing the majorants,
which are very small in comparison to the
computational times of the solver.

In {\bf Example 1}, we set the desired state
\begin{align*}
 y_d(\boldsymbol{x},t) = \frac{e^t \sin (t)}{10} \left(\left(12+4 \pi ^4\right) \sin ^2(t)-6 \cos ^2(t)-6 \sin (t) \cos (t)\right)
			 \sin (x_1 \pi) \sin (x_2 \pi),
\end{align*}
where $T=2 \pi / \omega$ and $\omega = 1$.
This function is time-periodic and analytic, but not time-harmonic.
Hence, the truncated Fourier series approximation of 
$y_d$ has to be computed for
applying the MhFEM as presented in Section \ref{Sec3:MhFEApprox}.
For that, the Fourier coefficients of $y_d$ can be computed analytically, 
and, then, they are approximated by the FEM. 
Next the systems (\ref{equation:MultiFESysBlockk}) and (\ref{equation:MultiFESysBlock0})
are solved for all $k \in \{0,\dots,N\}$. We choose the truncation index $N=8$. 
Finally, we reconstruct the fluxes by a $RT^0$-extension 
and compute the corresponding majorants.
The exact state is given by
$ y(\boldsymbol{x},t) = e^t \sin(t)^3 \sin (x_1 \pi) \sin (x_2 \pi) $.
In Table~\ref{tab:Ex1:0}, we present the CPU times $t^{\text{sec}}$, 
the majorants
  \begin{align*}
  \begin{aligned}
  \mathcal{M}^{\oplus_0}_{|\cdot|}
  = \sqrt{2} &\Big(C_F \, \big(\|{\mathcal{R}_1}^c_0
   \|_{\Omega}
   + \|{\mathcal{R}_3}^c_0
   \|_{\Omega} 
   \big)
   + \|{\mathcal{R}_2}^c_0
   \|_{\Omega}
    + \|{\mathcal{R}_4}^c_0
    \|_{\Omega} \Big) 
  \end{aligned}
 \end{align*}
 and $\mathcal{J}^{\oplus}_0$ as defined in (\ref{def:Joplus0}) as well as
the corresponding efficiency indices
\begin{align*}
 I_{\text{eff}}^{\mathcal{M},0} =
 \frac{\mathcal{M}^{\oplus_0}_{|\cdot|}}{|(y_0^c,p_0^c) - (\eta_0^c,\zeta_0^c)|_{1,\Omega}}
\end{align*}
and $I_{\text{eff}}^{\mathcal{J},0} = \mathcal{J}^{\oplus}_0 / \mathcal{J}_0$
obtained on grids of different mesh sizes.
\begin{table}[!ht]
\begin{center}
\begin{tabular}{|c|ccccc|}
  \hline
   grid & $t^{\text{sec}}$ 
   & $\mathcal{M}^{\oplus_0}_{|\cdot|}$ 
   & $I_{\text{eff}}^{\mathcal{M},0}$ 
   & $\mathcal{J}^{\oplus}_0$ 
   & $I_{\text{eff}}^{\mathcal{J},0}$ \\
  \hline
   $16   \times   16$   & 0.02   & 1.75e+01 & 2.50 & 1.26e+05 & 1.01 \\  
   $32   \times   32$   & 0.08   & 8.20e+00 & 2.20 & 1.27e+05 & 1.00 \\   
   $64   \times   64$   & 0.35   & 3.92e+00 & 2.05 & 1.27e+05 & 1.00 \\  
   $128 \times 128$   & 1.62   & 1.91e+00 & 1.98 & 1.27e+05 & 1.00  \\  
   $256 \times 256$   & 7.00   & 9.44e-01 & 1.94 & 1.27e+05 &  1.00 \\  
  \hline
\end{tabular}
\end{center}
\caption{The majorants $\mathcal{M}^{\oplus_0}_{|\cdot|}$ and $\mathcal{J}^{\oplus}_0$, and
their efficiency indices (Example 1).}
\label{tab:Ex1:0}
\end{table}
Moreover, Table~\ref{tab:Ex1:1} presents 
the CPU times $t^{\text{sec}}$, the majorants
 \begin{align*}
  \begin{aligned}
  \mathcal{M}^{\oplus_k}_{|\cdot|}
  = \sqrt{2} \Big(C_F \, 
   \big(\|{\mathcal{R}_1}_k \|_{\Omega}
      + \|{\mathcal{R}_3}_k \|_{\Omega}
	\big) +\, \|{\mathcal{R}_2}_k\|_{\Omega}
   + \|{\mathcal{R}_4}_k
   \|_{\Omega}  \Big)
  \end{aligned}
 \end{align*}
and $\mathcal{J}^{\oplus}_k$ as defined in (\ref{def:Joplusk})
for $k=1$ and, finally, the corresponding efficiency indices
\begin{align*}
 I_{\text{eff}}^{\mathcal{M},k} =
 \frac{\mathcal{M}^{\oplus_k}_{|\cdot|}}{|(\boldsymbol{y}_k,\boldsymbol{p}_k)-
(\boldsymbol{\eta}_k,\boldsymbol{\zeta}_k)|_{1,\Omega}}
\end{align*}
and $I_{\text{eff}}^{\mathcal{J},k} = \mathcal{J}^{\oplus}_k / \mathcal{J}_k$
obtained on grids of different mesh sizes.
Similar results are obtained for larger $k$ as well, which is illustrated
in Table~\ref{tab:Ex1:global}.
This table compares
the results for the modes $k \in \{0,\dots,8\}$
computed on the $256 \times 256$-mesh
and presents
the overall functional error estimates.
For that,
the remainder term
$\mathcal{E}_N$
is precomputed exactly,
see Remark \ref{remark:remainderterm}.
It can be observed that the values of the majorants 
$\mathcal{M}^{\oplus_k}_{|\cdot|}$ and
$\mathcal{J}^{\oplus}_{k}$
decrease for increasing $k$,
but that the values of the efficiency indices 
are all about the same,
which is a demonstration for 
the \textit{robustness} of the method with respect to the modes. 
Note that the overall efficiency index for $N=6$ is large ($I_{\text{eff}}^{\mathcal{M}} = 3.15$) 
compared to the efficiency indices corresponding to the single modes. The reason for that is
that the remainder term $\mathcal{E}_6 = 640.25$ is still quite large, and hence, more modes are needed.
For $N=8$, the remainder term is $\mathcal{E}_8 =106.07$,
which leads to a much better overall efficiency index ($I_{\text{eff}}^{\mathcal{M}} = 1.69$). 
The value for the cost functional is however sufficiently accurate with $N=6$.
This example demonstrates 
that the {\it a posteriori} error estimates clearly show what amount of
modes would be sufficient for representing the solution with a desired accuracy.
\begin{table}[!ht]
\begin{center}
\begin{tabular}{|c|ccccc|}
  \hline
   grid 
   & $t^{\text{sec}}$ 
   & $\mathcal{M}^{\oplus_1}_{|\cdot|}$ 
   & $I_{\text{eff}}^{\mathcal{M},1}$
   & $\mathcal{J}^{\oplus}_1$
   & $I_{\text{eff}}^{\mathcal{J},1}$  \\
  \hline
   $16   \times   16$   & 0.02   & 3.40e+01 & 2.50 & 4.74e+05 & 1.00  \\  
   $32   \times   32$   & 0.09   & 1.59e+01 & 2.20 & 4.79e+05 & 1.00 \\   
   $64   \times   64$   & 0.36   & 7.63e+00 & 2.05 & 4.80e+05 & 1.00 \\  
   $128 \times 128$   & 1.62   & 3.72e+00 & 1.98 & 4.80e+05 & 1.00  \\  
   $256 \times 256$   & 6.99   & 1.84e+00 & 1.94 & 4.80e+05 & 1.00 \\   
  \hline
\end{tabular}
\end{center}
\caption{The majorants $\mathcal{M}^{\oplus_1}_{|\cdot|}$ and $\mathcal{J}^{\oplus}_1$, and
their efficiency indices (Example 1).}
\label{tab:Ex1:1}
\end{table}
\begin{table}[!ht]
\begin{center}
\begin{tabular}{|c|ccccc|}
  \hline
   & $t^{\text{sec}}$ 
   & $\mathcal{M}^{\oplus}_{|\cdot|}$ 
   & $I_{\text{eff}}^{\mathcal{M}}$
   & $\mathcal{J}^{\oplus}$
   & $I_{\text{eff}}^{\mathcal{J}}$  \\
  \hline
   $k=0$   &   7.00 & 9.44e-01  & 1.94 & 1.27e+05 & 1.00 \\  
   $k=1$   &   6.99 & 1.84e+00 & 1.94 & 4.80e+05 & 1.00 \\    
   $k=2$   &   7.02 & 1.18e+00 & 1.94 & 1.99e+05 & 1.00 \\     
   $k=3$   &   7.17 & 6.78e-01 & 1.94 & 6.74e+04 & 1.00 \\    
   $k=4$   &   6.95 & 2.35e-01 & 1.92 & 8.42e+03 & 1.00 \\    
   $k=5$   &   7.01 & 8.57e-02 & 1.94 & 1.13e+03 & 1.00 \\    
   $k=6$   &   6.70 & 4.03e-02 & 1.87 & 2.29e+02 & 1.00 \\    
   $k=7$   &   6.77 & 2.13e-02 & 2.12 & 6.37e+01 & 1.03 \\       
   $k=8$   &   6.87 & 1.25e-02 & 2.18 & 2.19e+01 & 1.04 \\       
  \hline
   overall with $N=6$ & -- & 1.03e+01 & 3.15 & 3.17e+06 & 1.00 \\    
   overall with $N=8$ & -- & 6.13e+00 & 1.69 & 3.17e+06 & 1.00 \\     
  \hline
\end{tabular}
\end{center}
\caption{The overall majorants $\mathcal{M}^{\oplus}_{|\cdot|}$ and $\mathcal{J}^{\oplus}$,
and their parts computed on a $256 \times 256$-mesh (Example 1).}
\label{tab:Ex1:global}
\end{table}

In {\bf Example 2}, we set 
\begin{align*}
 y_d(\boldsymbol{x},t) = \frac{e^t}{10}  \left(-2 \cos(t)+ (10 + 4 \pi^4) \sin (t) \right)
			 \sin (x_1 \pi) \sin (x_2 \pi),
\end{align*}
where $T=2 \pi / \omega$ with $\omega = 1$. 
It is easy to see that this function is time-analytic, but not time-periodic.
As in the first example, we 
compute the MhFE approximation
of the desired state and solve the systems 
(\ref{equation:MultiFESysBlockk}) and (\ref{equation:MultiFESysBlock0})
for all $k \in \{0,\dots,N\}$ with first $N=6$ and second $N=8$ being the truncation index.
Finally, we compute also the solutions for $N=10$.
The exact state is given by $ y(\boldsymbol{x},t) = e^t \sin(t) \sin (x_1 \pi) \sin (x_2 \pi) $.
The results related to computational expenditures and efficiency indices are quite similar
to those for Example 1. Therefore, we present only numerical results in the form similar to
Table \ref{tab:Ex1:global} (see
Table~\ref{tab:Ex2:global}). 
In this numerical experiment, we again 
observe good and satisfying efficiency indices for $\mathcal{M}^{\oplus_k}_{|\cdot|}$. 
The remainder terms for $N=6$, $N=8$ and $N=10$ are 
$\mathcal{E}_6 = 44094.84$, $\mathcal{E}_8 = 19869.30$ and 
$\mathcal{E}_{10} = 10597.20$, respectively.
The efficiency index for $\mathcal{M}^{\oplus}_{|\cdot|}$ with $N=10$ improves a lot compared to
the index with $N=6$. 
Note that -- as in the first example -- the efficiency indices for the cost functional are approximately one.
This again demonstrates the accurateness of the majorants for the cost functional.

\begin{table}[!ht]
\begin{center}
\begin{tabular}{|c|ccccc|}
  \hline
   & $t^{\text{sec}}$ & $\mathcal{M}^{\oplus}_{|\cdot|}$ 
   & $I_{\text{eff}}^{\mathcal{M}}$
   & $\mathcal{J}^{\oplus}$
   & $I_{\text{eff}}^{\mathcal{J}}$  \\
  \hline
   $k=0$   &  6.86  & 1.58e+00 & 1.95 & 3.56e+05 & 1.00 \\  
   $k=1$   &  6.89  & 2.83e+00 & 1.95 & 1.14e+06 & 1.00 \\    
   $k=2$   &  6.79  & 1.41e+00 & 1.93 & 2.85e+05 & 1.00 \\     
   $k=3$   &  6.86  & 6.73e-01  & 1.89 & 6.69e+04 & 1.00 \\       
   $k=4$   &  6.76  & 3.78e-01  & 1.85 & 2.19e+04 & 1.00 \\      
   $k=5$   &  6.93  & 2.37e-01  & 1.78 & 9.05e+03 & 1.00 \\         
   $k=6$   &  6.83  & 1.60e-01  & 1.71 & 4.38e+03 & 1.00 \\        
   $k=7$   &  6.73  & 1.12e-01  & 1.62 & 2.37e+03 & 1.00 \\           
   $k=8$   &  6.99  & 8.28e-02  & 1.54 & 1.40e+03 & 1.00 \\         
   $k=9$   &  6.83  & 6.28e-02  & 1.45 & 8.74e+02 & 1.00 \\  
   $k=10$ &  6.87  & 4.88e-02  & 1.37 & 5.75e+02 & 1.00 \\     
  \hline
   overall with $N=6$   & -- & 6.94e+02 & 2.56 & 7.06e+06 & 1.00 \\      
   overall with $N=8$   & -- & 4.75e+02 & 1.75 & 7.06e+06 & 1.00 \\         
   overall with $N=10$ & -- & 3.55e+02 & 1.31 & 7.06e+06 & 1.00 \\    
  \hline
\end{tabular}
\end{center}
\caption{The overall majorants $\mathcal{M}^{\oplus}_{|\cdot|}$ and $\mathcal{J}^{\oplus}$,
and their parts computed on a $256 \times 256$-mesh (Example 2).}
\label{tab:Ex2:global}
\end{table}

In {\bf Example 3}, we set 
\begin{align*}
 y_d(\boldsymbol{x},t) = 
 \chi_{[\frac{1}{4},\frac{3}{4}]}(t) \,\chi_{[\frac{1}{2},1]^2}(\boldsymbol{x}),
\end{align*}
where $\chi$ denotes the characteristic function 
in space and time.
Let $T=1$, then $\omega = 2 \pi$.
Again the coefficients of the Fourier 
expansion associated with $y_d$ 
can be found analytically. They are
\begin{align*}
 {y_d^c}_k(\boldsymbol{x})  = \frac{
 \left(-\sin(\frac{k \pi}{2}) 
 + \sin(\frac{3 k \pi}{2}) \right)}{k \pi} \, \chi_{[\frac{1}{2},1]^2}(\boldsymbol{x}),
\end{align*}
and ${y_d^s}_k(\boldsymbol{x}) = 0$ for all $k \in \mathbb{N}$.
For $k=0$, ${y_d^c}_0(\boldsymbol{x}) = \chi_{[\frac{1}{2},1]^2}(\boldsymbol{x})/2$.
Since the exact solution cannot be 
computed analytically, we compute its MhFE approximation on a finer mesh 
($512 \times 512$-mesh). 
Since the modes ${y_d^c}_k(\boldsymbol{x}) = 0$ for all even $k \in \mathbb{N}$,
it suffices to show the results for odd modes as well as for $k=0$.
Table \ref{tab:Ex3:global} presents the results with 
truncation index $N=23$, since the results regarding the efficiency indices are similar
for higher modes. 
The computational times presented include 
the times for computing the
approximations of the exact modes on the finer mesh.
The numerical results again show 
the efficiency of the majorants for both, the discretization error and the cost functional.
This is especially observed for the majorant related to the cost functional, which is
very close to the exact value (in spite of a really complicated $y_d$).
The majorant $\mathcal{M}^{\oplus}_{|\cdot|}$ exposes an overestimation
but anyway provides realistic estimates of the errors in the state and control functions
measured in terms of the combined error norm.
\begin{table}[!ht]
\begin{center}
\begin{tabular}{|c|ccccc|}
  \hline
   & $t^{\text{sec}}$ 
   & $\mathcal{M}^{\oplus_k}_{|\cdot|}$ 
   & $I_{\text{eff}}^{\mathcal{M}}$
   & $\mathcal{J}^{\oplus_k}$
   & $I_{\text{eff}}^{\mathcal{J}}$  \\
  \hline
   $k=0$   &  38.60 
     & 3.96e+02 & 3.64 & 8.20e+04 & 1.32 \\ 
   $k=1$   &  38.88 
     & 4.80e+02 & 3.73 & 1.31e+05 & 1.30 \\    
   $k=3$   &  38.82 
     & 1.22e+02 & 2.52 & 1.35e+04 & 1.21 \\        
   $k=5$   &  38.98 
     & 5.58e+01 & 2.41 & 4.62e+03 & 1.14 \\       
   $k=7$   &  38.64 
     & 3.22e+00 & 2.45 & 2.28e+03 & 1.11 \\           
   $k=9$   &  39.06 
     & 2.12e+01 & 2.51 & 1.35e+03 & 1.08 \\  
   $k=11$ &  38.92 
     & 1.50e+01 & 2.55 & 8.89e+02 & 1.07 \\   
   $k=13$ &  39.13 
     & 1.14e+01 & 2.62 & 6.30e+02 & 1.06 \\     
   $k=15$ &  38.59 
     & 8.93e+00 & 2.63 & 4.68e+02 & 1.04 \\     
   $k=17$ &  38.58 
     & 7.36e+00 & 2.70 & 3.63e+02 & 1.04 \\     
   $k=19$ &  38.78 
     & 6.06e+00 & 2.69 & 2.88e+02 & 1.03 \\     
   $k=21$ &  38.72 
     & 5.27e+00 & 2.78 & 2.36e+02 & 1.03 \\     
   $k=23$ &  38.87 
     & 4.47e+00 & 2.74 & 1.96e+02 & 1.02 \\     
  \hline
\end{tabular}
\end{center}
\caption{The majorants and corresponding efficiency indices 
computed on a $256 \times 256$-mesh (Example 3).}
\label{tab:Ex3:global}
\end{table}

\section{Conclusions and Outlook}
\label{Sec7:ConclusionsOutlook}

In \cite{LRW:LangerRepinWolfmayr:2015}, the authors derived functional-type 
{\it a posteriori} error estimates for MhFE 
approximations to linear parabolic time-periodic boundary value problems.
In this work, we extend this technique
to the derivation of {\it a posteriori} error estimates for MhFE  
solutions of the corresponding distributed optimal control problem
which leads to additional challenges in the analysis.
The reduced optimality system is nothing but a coupled 
parabolic time-periodic PDE system for the state and 
the adjoint state. We are not only interested in computable
{\it a posteriori} error bounds for the state, the adjoint state 
and the control, but also for the cost functional.
In case of linear time-periodic parabolic constraints, 
the approximation via MhFE functions 
leads to the decoupling of computations related to different modes.
Due to this feature of the MhFEM, 
we can in principle use different meshes for different modes and 
independently generate them by adaptive
finite element approximations of the respective Fourier coefficients.
To assure the quality of approximations constructed in this way,
we need
fully computable {\it a posteriori} estimates, 
which provide guaranteed bounds of global 
errors and reliable indicators of errors associated with the modes.
Then, by prescribing certain bounds, we can finally
filter out the Fourier coefficients, which are important for the
numerical solution of the problem.
This technology will lead to an 
\emph{adaptive multiharmonic finite element method (AMhFEM)}
that will provide complete adaptivity in space and time.
The development and the analysis of such an  AMhFEM 
goes beyond the scope of this paper,
but will heavily be based on the results of this paper 
as described above. 
It is clear that the functional {\it a posteriori} estimates derived here for
time-harmonic parabolic optimal control problems can also be
obtained for distributed time-harmonic eddy current  optimal control problems
as studied in 
\cite{PhD:Kolmbauer:2012c, PhD:KolmbauerLanger:2012, PhD:KolmbauerLanger:2013}.

\section*{Acknowledgments}

The authors gratefully acknowledge the financial support by the Austrian Science Fund 
(FWF) under the grants NST-0001 and W1214-N15 (project DK4), the Johannes Kepler University of Linz and
the Federal State of Upper Austria.

\bibliographystyle{abbrv}
\bibliography{LRW2arxiv.bib}

\begin{thebibliography}{10}

\bibitem{LRW:Altmann:2013}
K.~Altmann.
\newblock {\em Numerische Verfahren der Optimalsteuerung von Magnetfeldern}.
\newblock PhD thesis, Technical University of Berlin, 2013.

\bibitem{LRW:AltmannStingelinTroeltzsch:2014}
K.~Altmann, S.~Stingelin, and F.~Tr{\"o}ltzsch.
\newblock On some optimal control problems for electrical circuits.
\newblock {\em International Journal of Circuit Theory and Applications},
  42(8):808--830, 2014.

\bibitem{LRW:BachingerKaltenbacherReitzinger:2002a}
F.~Bachinger, M.~Kaltenbacher, and S.~Reitzinger.
\newblock An efficient solution strategy for the {HBFE} method.
\newblock {\em Proceedings of the IGTE '02 Symposium Graz, Austria}, pages
  385--389, 2002.

\bibitem{PhD:BachingerLangerSchoeberl:2005}
F.~Bachinger, U.~Langer, and J.~Sch\"oberl.
\newblock Numerical analysis of nonlinear multiharmonic eddy current problems.
\newblock {\em Numer. Math.}, 100(4):593--616, 2005.

\bibitem{PhD:BachingerLangerSchoeberl:2006}
F.~Bachinger, U.~Langer, and J.~Sch\"oberl.
\newblock Efficient solvers for nonlinear time-periodic eddy current problems.
\newblock {\em Comput. Vis. Sci.}, 9(4):197--207, 2006.

\bibitem{LRW:BorziSchulz:2012}
A.~Borz\`{i} and V.~Schulz.
\newblock {\em Computational Optimization of Systems Governed by Partial
  Differential Equations}.
\newblock SIAM, Philadelphia, 2012.

\bibitem{LRW:Braess:2005}
D.~Braess.
\newblock {\em Finite elements: Theory, fast solvers, and applications in solid
  mechanics}.
\newblock Cambridge University Press, second edition, 2005.

\bibitem{LRW:BrezziFortin:1991}
F.~Brezzi and M.~Fortin.
\newblock {\em Mixed and hybrid finite element methods}, volume~15 of {\em
  Springer Series in Computational Mathematics}.
\newblock Springer, New York, 1991.

\bibitem{LRW:Ciarlet:1978}
P.~G. Ciarlet.
\newblock {\em The Finite Element Method for Elliptic Problems}, volume~4 of
  {\em Studies in Mathematics and its Applications}.
\newblock North-Holland, Amsterdam, 1978.
\newblock Republished by SIAM in 2002.

\bibitem{PhD:CopelandLanger:2010}
D.~M. Copeland and U.~Langer.
\newblock Domain decomposition solvers for nonlinear multiharmonic finite
  element equations.
\newblock {\em J. Numer. Math.}, 18(3):157--175, 2010.

\bibitem{PhD:GaevskayaHoppeRepin:2006}
A.~Gaevskaya, R.~H.~W. Hoppe, and S.~Repin.
\newblock A posteriori estimates for cost functionals of optimal control
  problems.
\newblock {\em Numerical Mathematics and Advanced Applications, Proceedings of
  the ENUMATH 2005}, pages 308--316, 2006.

\bibitem{PhD:GaevskayaHoppeRepin:2007}
A.~Gaevskaya, R.~H.~W. Hoppe, and S.~Repin.
\newblock Functional approach to a posteriori error estimation for elliptic
  optimal control problems with distributed control.
\newblock {\em J. Math. Sci.}, 144(6):4535--4547, 2007.

\bibitem{PhD:GaevskayaRepin:2005}
A.~V. Gaevskaya and S.~I. Repin.
\newblock A posteriori error estimates for approximate solutions of linear
  parabolic problems.
\newblock {\em Differential Equations}, 41(7):970--983, 2005.

\bibitem{LRW:HinzePinnauUlbrichUlbrich:2009}
M.~Hinze, R.~Pinnau, M.~Ulbrich, and S.~Ulbrich.
\newblock {\em Optimization with PDE constraints}.
\newblock Mathematical Modelling: Theory and Applications 23, Springer, Berlin,
  2009.

\bibitem{LRW:HouskaLogistImpeDiehl:2009}
B.~Houska, F.~Logist, J.~V. Impe, and M.~Diehl.
\newblock Approximate robust optimization of time-periodic stationary states
  with application to biochemical processes.
\newblock In {\em Proceedings of the 48th IEEE Conference on Decision and
  Control}, pages 6280--6285. Shanghai, China, 2009.

\bibitem{LRW:JungLanger:2013}
M.~Jung and U.~Langer.
\newblock {\em Methode der finiten Elemente f\"ur Ingenieure: Eine Einf\"uhrung
  in die numerischen Grundlagen und Computersimulation}.
\newblock Springer, Wiesbaden, second edition, 2013.

\bibitem{PhD:KollmannKolmbauer:2011}
M.~Kollmann and M.~Kolmbauer.
\newblock A preconditioned {M}in{R}es solver for time-periodic parabolic
  optimal control problems.
\newblock {\em Numer. Linear Algebra Appl.}, 20(5):761--784, 2013.

\bibitem{LRW:KollmannKolmbauerLangerWolfmayrZulehner:2013}
M.~Kollmann, M.~Kolmbauer, U.~Langer, M.~Wolfmayr, and W.~Zulehner.
\newblock A finite element solver for a multiharmonic parabolic optimal control
  problem.
\newblock {\em Comput. Math. Appl.}, 65(3):469--486, 2013.

\bibitem{PhD:Kolmbauer:2012c}
M.~Kolmbauer.
\newblock {\em The Multiharmonic Finite Element and Boundary Element Method for
  Simulation and Control of Eddy Current Problems}.
\newblock PhD thesis, JKU Linz, 2012.

\bibitem{PhD:KolmbauerLanger:2012}
M.~Kolmbauer and U.~Langer.
\newblock A robust preconditioned {M}in{R}es solver for distributed
  time-periodic eddy current optimal control problems.
\newblock {\em SIAM J. Sci. Comput.}, 34(6):B785--B809, 2012.

\bibitem{PhD:KolmbauerLanger:2013}
M.~Kolmbauer and U.~Langer.
\newblock Efficient solvers for some classes of time-periodic eddy current
  optimal control problems.
\newblock {\em Numerical Solution of Partial Differential Equations: Theory,
  Algorithms, and Their Applications}, 45:203--216, 2013.

\bibitem{LRW:Kraus:2012}
J.~Kraus.
\newblock Additive {S}chur complement approximation and application to
  multilevel preconditioning.
\newblock {\em SIAM J. Sci. Comput.}, 34(6):A2872--A2895, 2012.

\bibitem{LRW:KrausWolfmayr:2013}
J.~Kraus and M.~Wolfmayr.
\newblock On the robustness and optimality of algebraic multilevel methods for
  reaction-diffusion type problems.
\newblock {\em Comput. Vis. Sci.}, 16(1):15--32, 2013.

\bibitem{PhD:KrendlSimonciniZulehner:2012}
W.~Krendl, V.~Simoncini, and W.~Zulehner.
\newblock Stability estimates and structural spectral properties of saddle
  point problems.
\newblock {\em Numer. Math.}, 124(1):183--213, 2013.

\bibitem{LRW:Ladyzhenskaya:1973}
O.~A. Ladyzhenskaya.
\newblock {\em The Boundary Value Problems of Mathematical Physics}.
\newblock Nauka, Moscow, 1973.
\newblock In Russian. Translated in \emph{Appl. Math. Sci.} 49, Springer, 1985.

\bibitem{LRW:LadyzhenskayaSolonnikovUralceva:1968}
O.~A. Ladyzhenskaya, V.~A. Solonnikov, and N.~N. Ural'ceva.
\newblock {\em Linear and Quasilinear Equations of Parabolic Type}.
\newblock AMS, Providence, RI, 1968.

\bibitem{LRW:LangerRepinWolfmayr:2015}
U.~Langer, S.~Repin, and M.~Wolfmayr.
\newblock Functional a posteriori error estimates for parabolic time-periodic
  boundary value problems.
\newblock {\em Comput. Methods Appl. Math.}, 15(3):353--372, 2015.

\bibitem{LRW:LangerWolfmayr:2013}
U.~Langer and M.~Wolfmayr.
\newblock Multiharmonic finite element analysis of a time-periodic parabolic
  optimal control problem.
\newblock {\em J. Numer. Math.}, 21(4):265--300, 2013.

\bibitem{PhD:MaliNeittaanmaekiRepin:2014}
O.~Mali, P.~Neittaanm\"{a}ki, and S.~Repin.
\newblock {\em Accuracy Verification Methods. Theory and Algorithms}, volume~32
  of {\em Computational Methods in Applied Sciences}.
\newblock Springer, Netherlands, 2014.

\bibitem{LRW:NeittaanmaekiSprekelsTiba:2006}
P.~Neittaanm\"{a}ki, J.~Sprekels, and D.~Tiba.
\newblock {\em Optimization of Elliptic Systems: Theory and Applications}.
\newblock Springer, 2006.

\bibitem{PhD:RaviartThomas:1977}
P.~A. Raviart and J.~M. Thomas.
\newblock A mixed finite element method for 2-nd order elliptic problems.
\newblock {\em Mathematical Aspects of Finite Element Methods, Lect. Notes
  Math.}, 606:292--315, 1977.

\bibitem{PhD:Repin:2002}
S.~Repin.
\newblock Estimates of deviation from exact solutions of initial-boundary value
  problems for the heat equation.
\newblock {\em Rend. Mat. Acc. Lincei}, 13(2):121--133, 2002.

\bibitem{PhD:Repin:2008}
S.~Repin.
\newblock {\em A Posteriori Estimates for Partial Differential Equations}.
\newblock Radon Series on Computational and Applied Mathematics 4, Walter de
  Gruyter, Berlin, 2008.

\bibitem{LRW:RobertsThomas:1991}
J.~E. Roberts and J.~M. Thomas.
\newblock Mixed and hybrid methods.
\newblock {\em Handbook of numerical analysis}, 2(1):523--639, 1991.

\bibitem{LRW:Steinbach:2008}
O.~Steinbach.
\newblock {\em Numerical Approximation Methods for Elliptic Boundary Value
  Problems: Finite and Boundary Elements}.
\newblock Springer, New York, 2008.

\bibitem{LRW:Troeltzsch:2010}
F.~Tr\"oltzsch.
\newblock {\em Optimal Control of Partial Differential Equations. Theory,
  Methods and Applications}.
\newblock Graduate Studies in Mathematics 112, AMS, Providence, RI, 2010.

\bibitem{LRW:Wolfmayr:2014}
M.~Wolfmayr.
\newblock {\em Multiharmonic Finite Element Analysis of Parabolic Time-Periodic
  Simulation and Optimal Control Problems}.
\newblock PhD thesis, Johannes Kepler University of Linz, 2014.

\bibitem{PhD:YamadaBessho:1988}
S.~Yamada and K.~Bessho.
\newblock Harmonic field calculation by the combination of finite element
  analysis and harmonic balance method.
\newblock {\em IEEE Trans. Magn.}, 24(6):2588--2590, 1988.

\end{thebibliography}
\end{document}